\documentclass[10pt,a4paper]{amsart}

 \usepackage{amssymb}   
 \usepackage{amstext}
 \usepackage{amsfonts}
 \usepackage{amsmath}   
 \usepackage{amscd}     
 \usepackage{latexsym}  
 \usepackage{color}
 \usepackage{epsfig}
 \usepackage{epic}
 \usepackage{pstricks,pst-node,pst-text,pst-tree}
 \usepackage[all]{xy}
 \usepackage{array}
 \usepackage{amsthm}
 \usepackage{eucal, mathrsfs} 
 
\newcommand{\numberset}[1]{\ensuremath{\mathbb{#1}}}    
\newcommand{\C}{\numberset{C}}  
\newcommand{\N}{\numberset{N}}  
\newcommand{\R}{\numberset{R}}  
\newcommand{\Z}{\numberset{Z}}  
\newcommand{\PP}{\numberset{P}}  

\newcommand{\zbar}{\ensuremath{\overline{z}}}

\newcommand{\inner}[2]{ \langle {#1}, {#2} \rangle}
\newcommand{\znor}{Z_{\text{nor}}}
\newcommand{\zbnor}{\bar{Z}_{\text{nor}}}
\newcommand{\gnor}{\Gamma_{\text{nor}}}

\theoremstyle{definition}

\newtheorem{thm}{Theorem}[section]
\newtheorem{prop}[thm]{Proposition}
\newtheorem{lem}[thm]{Lemma}
\newtheorem{cor}[thm]{Corollary}
\newtheorem{rem}[thm]{Remark}
\newtheorem{ex}[thm]{Example}
\newtheorem{defi}[thm]{Definition}

\newtheorem{ass}[thm]{Assumption}

\DeclareMathOperator{\re}{Re}

\DeclareMathOperator{\Crit}{Crit}

\DeclareMathOperator{\Gl}{Gl}
\DeclareMathOperator{\spn}{span}

\DeclareMathOperator{\I}{Id}
\DeclareMathOperator{\inv}{inv}
\DeclareMathOperator{\Log}{Log}
\addtocounter{secnumdepth}{-2}

\author{R. Casta\~no-Bernard and D. Matessi}
\title{Semi-global invariants of piecewise smooth Lagrangian fibrations}

\begin{document}
\maketitle

\begin{abstract}
We study certain types of piecewise smooth Lagrangian fibrations of smooth symplectic manifolds, which we call \emph{stitched Lagrangian fibrations}. We extend the classical theory of action-angle coordinates to these fibrations by encoding the information on the  non-smoothness into certain invariants consisting, roughly, of a sequence of closed $1$-forms on a torus. The main motivation for this work is given by the piecewise smooth Lagrangian fibrations previously constructed by the authors \cite{CB-M-torino}, which topologically coincide with the local models used by Gross in Topological Mirror Symmetry \cite{TMS}. 
\end{abstract}

\section{Introduction}

Lagrangian fibrations arise naturally from integrable systems. It is a
standard fact of Hamiltonian mechanics that such fibrations are locally given by
maps of the type:
\[
f=(f_1,\ldots, f_n),
\]
where the function components of $f$ are Poisson commuting functions on a
symplectic manifold and such that the differentials $df_1,\ldots ,df_n$ are
pointwise linearly independent almost everywhere. 
It is customary to assume $f$ to be $C^\infty$ differentiable (smooth). Under this regularity assumption, a
classical theorem of Arnold-Liouville says that a smooth proper
Lagrangian submersion with connected fibres has locally the structure of a
trivial Lagrangian $T^n$-bundle. In particular, all proper Lagrangian
submersions are locally modelled on $U\times T^n$, where $U\subseteq\R^n$ is a
contractible open set and $U\times T^n$ has the standard symplectic form induced 
from $\R^{2n}$. Standard coordinates with values in $U\times T^n$
are known as action-angle coordinates. Since these are defined on a fibred
neighbourhood, action-angle coordinates are \emph{semi-global} canonical
coordinates. Thus proper Lagrangian submersions have no semi-global
symplectic invariants. 

\medskip
In this article we investigate the semi-global symplectic topology of
proper Lagrangian fibrations given by piecewise smooth maps.   In \cite{CB-M-torino}\S 6 we
introduced the notion of \textit{stitched Lagrangian fibration}. These are continuous
proper $S^1$ invariant fibrations of smooth symplectic manifolds $X$ which fail to
be smooth only along the zero level set $Z=\mu^{-1}(0)$ of the moment map of the $S^1$
action and whose fibres are all smooth Lagrangian $n$-tori. Essentially, these
fibrations consist of two honest smooth pieces $X^+=\{\mu\geq 0\}$ and $X^{-}=\{\mu\leq 0\}$, stitched\footnote{We
have chosen to use `stitching' rather than `gluing' since the resulting map is
in general non smooth; the term `gluing' usually has a smoothness meaning
attached to it.} together along  $Z$,  which we call the \textit{seam}. These fibrations, roughly speaking, can be expressed locally as:
\[
f=(\mu, f_2^\pm,\ldots, f_n^\pm),
\]
where $f_j^+$ and $f_j^-$ are smooth functions defined on $X^+$ and $X^-$, respectively, whose differentials do not necessarily coincide along $Z$. Fibrations of this type are
implicit in the examples proposed earlier by the authors \cite{CB-M-torino}\S
5 and may also be implicit in those in \cite{Ruan}. In this paper, we develop a theory of  action-angle coordinates for this class of piecewise smooth fibrations. Contrary to what happens in the smooth case, we found that these fibrations do give rise to semi-global symplectic invariants.

\medskip
To the authors' knowledge, the kind of non-smoothness we investigate here does not seem to be of relevance to Hamiltonian mechanics. Nevertheless it is an important issue in symplectic topology and mirror symmetry. Over the past ten years, Lagrangian torus fibrations, in particular those which are \emph{special Lagrangian}, have been discovered to play a fundamental role in mirror symmetry \cite{SYZ}.  One should expect mirror pairs of Calabi-Yau manifolds to be fibred by Lagrangian tori and the mirror relation to be expressed in terms of a Legendre transform between the corresponding affine bases \cite{Hitchin}, \cite{G-Siebert}, \cite{Kontsevich-Soibelman}. This approach to mirror symmetry has some intricacies. For instance, there are examples of (non proper) special Lagrangian fibrations which are not given by smooth maps. Actually, one should expect a generic special Lagrangian fibration to be piecewise smooth \cite{Joyce-SYZ}. Non-smoothness may also arise even in the purely Lagrangian case. In fact, there are examples of Lagrangian torus fibrations of Calabi-Yau manifolds which are piecewise smooth \cite{Ruan}. This lack of regularity has two important consequences. In first place, the discriminant locus of the fibration --i.e. the set of points in the base corresponding to singular fibres-- may have codimension less than 2. Secondly, the base of the fibration may no longer carry the structure of an integral affine manifold away from the discriminant locus. In fact, the affine structure may break off not only along the discriminant --as it normally occurs in the smooth case-- but also along a larger set containing the discriminant. Under these circumstances, it may become problematic to interpret the SYZ duality as a Legendre transform between affine manifolds. One should therefore understand the symplectic topology of piecewise smooth Lagrangian fibrations. 

\medskip
Some of the piecewise smooth examples here actually resemble the singular behaviour expected to appear in generic special Lagrangian fibrations. What is more important for our purposes, however, is the fact that our Lagrangian models coincide topologically with the non-Lagrangian models used by Gross \cite{TMS}; the discriminant locus in our case may jump to codimension 1 in some regions but the total spaces are the same. In some cases, the discriminant has the shape of a planar amoeba $\Delta$ (see Figure~\ref{fig: amoeba}) and fails to be smooth over the hyperplane $\Gamma=\{\mu =0\}$ containing $\Delta$. Away from $\Delta$ these fibrations are stitched Lagrangian torus fibrations. In particular, the affine structure on the base breaks apart along $\Gamma\setminus\Delta$. In this paper we provide some useful techniques to understand how this degeneration of the affine structure occurs.

\medskip
The material of this paper is organised as follows. In \S \ref{section aa}
we start reviewing the classical theory of action-angle coordinates for
smooth fibrations. In \S \ref{sec:def-ex} we recall the construction of
piecewise smooth fibrations of \cite{CB-M-torino}, these are explicitly
given examples, some of them with codimension 1 discriminant locus. Then
we revise the definition of stitched fibration, introduced in
\cite{CB-M-torino}. We formalise the idea of action-angle coordinates for
stitched fibrations, allowing us to define the first order invariant, $\ell_1$, of a
stitched fibration. This invariant measures the discrepancy along $Z$ between the
distributions spanned by the Hamiltonian vector fields $\eta_2^+,\ldots ,\eta^+_n$ and $\eta_2^-,\ldots ,\eta^-_n$ corresponding to $f_2^+,\ldots ,f_n^+$ and $f_2^-,\ldots ,f_n^-$, respectively. 
The seam $Z$ is an $S^1$-bundle $p:Z\rightarrow \bar Z:=Z\slash S^1$ such that: 
\[
\xymatrix{
Z \ar[dr]_{f|_Z} \ar[rr]^{p} &  & \bar{Z} \ar[dl]^{\bar{f}} \\
&  \Gamma
}
\]
where $\bar f$ is the reduced fibration over the wall $\Gamma=\{\mu=0\}$, with tangent $(n-1)$-plane distribution: 
\[
\mathfrak{L}=\ker \bar f_\ast\subset T\bar Z.   
\]
Let $\mathscr L_{\bar Z}$ be the set of fibrewise closed sections of $\mathfrak{L}^\ast$, i.e. elements in $\mathscr L_{\bar Z}$ can be viewed as closed 1-forms on the fibres of $\bar f$. The first order invariant of $f$ is defined as follows. There are smooth $S^1$-invariant functions $a_2,\ldots ,a_n$ on $Z$ such that $a_j\eta_1=\eta_j^+-\eta_j^-$. In particular, this implies that $\eta_j^+$ and $\eta_j^-$ are mapped under $p_\ast$ to the same vector field $\bar\eta_j$ on $\bar Z$. The first order invariant $\ell_1$ is defined to be the section of $\mathfrak L^\ast$ such that $\ell_1(\bar\eta_j)=a_j$. It turns out that $\ell_1\in\mathscr L_{\bar Z}$.
In \S \ref{sec:ho-inv} we investigate higher order invariants. These are sequences $\{ \ell_k\}_{k\in\N}$, with each $\ell_k\in\mathscr L_{\bar Z}$. Given a stitched Lagrangian fibration $f$, we define $\inv (f)$ to consist of the data of $(\bar Z, \bar f)$, suitably normalized, together with the sequence $\ell=\{\ell_k\}_{k\in\N}$. The main result of this paper is proved in \S \ref{sec:normal} (cf. Theorem \ref{thm: grosso} and Theorem \ref{broken:constr2}) where we give a classification of stitched fibrations up to fibre-preserving symplectomorphims. Roughly, this can be stated as follows: 

\medskip
\noindent\textbf{Theorem.} There are stitched Lagrangian fibrations $f$ having any specified set of data $\inv (f)$. Moreover, given stitched Lagrangian fibrations $f$ and $f'$ with invariants $\inv(f)$ and $\inv(f')$, respectively, there is a smooth symplectomorphism $\Phi$, defined on a neighbourhood of $Z$, and a smooth diffeomorphism $\phi$ preserving $\Gamma\subset B$ and a commutative diagram:
\[
\begin{CD}
X @>\Phi>> X'\\
@VfVV @Vf'VV\\
B @>\phi>> B'
\end{CD}
\]
if and only if $\inv(f)=\inv(f')$.

\medskip
This result extends Arnold-Liouville's theorem to this piecewise smooth setting. In \S \ref{sec. stitch w/mono} we study stitched fibrations over non simply connected bases. We show that one can read the monodromy of a stitched fibration as a jump of the cohomology class $[\ell_1(b)]$ as $b\in\Gamma$ traverses a component of the discriminant locus. 

\medskip
In the last section, we propose the following:

\medskip
\noindent\textbf{Conjecture.} Let $Y\subseteq (\C^\ast)^{n-1}$ be a smooth algebraic hypersurface, $\Log :(\C^\ast)^{n-1}\rightarrow\R^{n-1}$ be the map defined by:
\[
\Log (z_2,\ldots ,z_n)=(\log |z_2|,\ldots ,\log |z_n|).
\] 
Then there is a piecewise smooth Lagrangian $n$-torus fibration with discriminant locus being the amoeba $\Delta=\Log (Y)$ inside $\{0\}\times\R^{n-1}\subset\R^n$. Away from $\Delta$ these fibrations are stitched Lagrangian fibrations.

\medskip
To support this conjecture we propose a construction.

\medskip
The results of this article allow us to have good control on the
regularity of a large class of proper Lagrangian fibrations. Using simple
techniques, one may deform the invariants of a given stitched fibration
and produce proper Lagrangian fibrations with $S^1$ symmetry which are
smooth on prescribed regions. This can be done, for instance, by multiplying a given sequence of invariants by a smooth function on the base $B$ vanishing on a prescribed region. In joint work in progress \cite{CB-M}, the
authors use these and other techniques to give a construction of
Lagrangian 3-torus fibrations of compact symplectic 6-manifolds starting from the information encoded in
suitable integral affine manifolds, such as those arising from toric degenerations \cite{G-Siebert2003}. Such affine structures are expected to appear as Gromov-Hausdorff limits of degenerating families of Calabi-Yau manifolds (in the sense of \cite{G-Wilson2}, \cite{Kontsevich-Soibelman}).

\section{Action-angle coordinates}\label{section aa}
We review the classical theory of action-angle coordinates for $C^\infty$ Lagrangian fibrations. For the details we refer the reader to \cite{Arnold}. Assume we are given a $2n$-dimensional symplectic manifold $X$ with symplectic structure $\omega$, a smooth $n$-dimensional manifold $B$ and a proper submersion $f:X\rightarrow B$ whose fibres are connected Lagrangian submanifolds. 

\medskip
Let $F_b$ be the fibre of $f$ over $b\in B$. We can define an action of $T_b^\ast B$ on $F_b$ as follows. For every $\alpha\in T^\ast_b B$ we can associate a vector field $v_\alpha$ on $F_b$ determined by
\begin{equation}\label{eq contr.}
\iota_{v_\alpha}\omega=f^\ast\alpha .
\end{equation}
Let $\phi_\alpha^t$ be the flow of $v_\alpha$ with time $t\in\R$. Define $\theta_\alpha$ as 
$\theta_\alpha(p)=\phi^1_\alpha (p)$ where $p\in F_b$. One can check that $\theta_\alpha$ is well defined and that it induces an action $(\alpha, p)\mapsto \theta_\alpha (p)$. Furthermore, the action is transitive. 
Then, $\Lambda_b$ defined as
\[ \Lambda_{b} = \{ \lambda \in T^{\ast}_{b}B \ | \ \theta_{\lambda}(p) = p,
\ \text{for all} \ p \in F_b \} \]
is a closed discrete subgroup of $T^\ast_bB$, i.e. a lattice.
From the properness of $f$ it follows that $\Lambda_b$ is maximal
(in particular homomorphic to $\Z^{n}$) and that $F_b$ is diffeomorphic
to $T^{\ast}_bB/ \Lambda_b$ and therefore $F_b$ is an $n$-torus.

\medskip
Let $\Lambda = \cup_{b \in B} \Lambda_b$. One can compute $\Lambda$ as follows.
Given a point $b_0 \in B$ and a contractible neighbourhood $U$ of $b_0$,
for every $b \in U$, $H_1(F_b, \Z)$ is naturally identified with 
$H_1(F_{b_0}, \Z)$. Choose a basis $\gamma_1, \ldots, \gamma_n$
of $H_1(F_{b_0},\Z)$. Given a vector field $v$ on 
$U$, denote by $\tilde{v}$ a lift of $v$ on $f^{-1}(U)$.
We can define the following $1$-forms $\lambda_1, \ldots, \lambda_n$
on $B$:
\begin{equation} \label{per:def}
 \lambda_j(v) = - \int_{\gamma_j} \iota_{\tilde{v}} \omega. 
\end{equation}
It is well known that the 1-forms $\lambda_j$ are closed and they generate $\Lambda$. If $\sigma :B\rightarrow X$ is a smooth section of $f$ we can define the map
\[
\Theta: T^\ast B\slash\Lambda\rightarrow X
\]
by $\Theta (b,\alpha )=\theta_\alpha(\sigma (b))$. This map is a diffeomorphism and  it is a symplectomorphism if $\sigma (B)\subseteq X$ is Lagrangian. A choice of functions $a_j$ such that $da_j=\lambda_j$ defines coordinates $a=(a_1,\ldots a_n)$ on $U$ called \textit{action coordinates}. In particular, a covering $\{U_i\}$ of $B$ by small enough contractible open sets and a choice of action coordinates $a_i$ on each $U_i$ defines an integral affine structure on $B$, i.e. an atlas whose change of coordinates maps are transformations in $\R^n\rtimes\Gl (n,\Z)$.

\medskip
A less invariant approach --but useful for explicit computations--  can be described as follows. Let $(b_1,\ldots b_n)$ be local coordinates on $U\subseteq B$ and let $f_j=b_j\circ f$. Then $f_1,\ldots ,f_n$ define an integrable Hamiltonian system. Let $\Phi^t_{\eta_j}$ be the flow of the Hamiltonian vector field $\eta_j$ of $f_j$. Let $\sigma$ be a Lagrangian section of $f$ over $U$. Then the map $\Theta$ above can be expressed as:
\[
\Theta: (b,t_1db_1+\cdots +t_ndb_n)\mapsto \Phi^{t_1}_{\eta_1}\circ\cdots\Phi^{t_n}_{\eta_n}(\sigma (b)).
\]
One may verify that
\[
\Lambda_b=\{ (b,t_1db_1+\cdots +t_ndb_n)\in T^\ast_bU\mid \Phi^{t_1}_{\eta_1}\circ\cdots\Phi^{t_n}_{\eta_n}(\sigma (b))=\sigma (b)\}.
\]
When $(b_1,\ldots , b_n)$ are action coordinates, $(b_1,\ldots , b_n, t_1,\ldots ,t_n)$ are \textit{action-angle coordinates}. These coordinates always exist on a fibred neighbourhood $f^{-1}(U)$ of a fibre $F_b$ with $U\subseteq\R^n$ a small neighbourhood of $b$, thus they can be regarded as \emph{semi-global} canonical coordinates. In particular, we have the following classical result:
\begin{thm}[Arnold, Liouville] A proper Lagrangian submersion with connected fibres and a Lagrangian section has no semi-global symplectic invariants.
\end{thm}
The \emph{global} existence of action-angle coordinates is obstructed. For the details concerning this issue we refer the reader to Duistermaat \cite{Dui}.

\medskip
In the next section we consider a larger class of Lagrangian submersions which include some Lagrangian fibrations which fail to be given by $C^\infty$ maps.
 
\section{Stitched Lagrangian fibrations: definitions and examples} 
                                                     \label{sec:def-ex}
 \begin{defi}\label{defi stitched} Let $(X, \omega)$ be a smooth $2n$-dimensional 
symplectic manifold. Suppose there is a free Hamiltonian $S^1$ action
on $X$ with moment map $\mu: X \rightarrow \R$. Let 
$X^+ = \{ \mu \geq 0 \}$ and $X^- = \{ \mu \leq 0 \}$. Given a smooth $(n-1)$-dimensional
manifold $M$, a map $f: X \rightarrow \R \times M$
is said to be a \textbf{stitched Lagrangian fibration} if 
there is a continuous $S^1$ invariant function $G: X \rightarrow M$, 
such that the following holds: 

\begin{itemize}
\item[\textit{(i)}] Let $G^{\pm} = G|_{X^{\pm}}$. Then $G^+$ and $G^-$ are restrictions of $C^\infty$ maps on $X$;
\item[\textit{(ii)}] $f$ can be written as \[ f =  (\mu, G)  \] and $f$ restricted to $X^{\pm}$ is a proper submersion with connected Lagrangian fibres. We denote 
\[ f^{\pm}=f|_{X^{\pm}}.\]
\end{itemize}
We call $Z = \mu^{-1}(0)$ the \textbf{seam}.
\end{defi}

We warn the reader that throughout the paper the superscript $\pm$ appearing in a 
sentence means that the sentence is true if read separately with the $+$ superscript 
and with the $-$ superscript. Notice that a stitched Lagrangian fibration may be non-smooth. In general it will be only piecewise $C^\infty$, however all its fibres are smooth Lagrangian tori. Observe also that $f^+$ and $f^-$ are restrictions of $C^\infty$ maps, they are not a priori required to extend to smooth Lagrangian fibrations beyond $X^+$ and $X^-$, respectively. Later we show, however, that for any stitched fibration, $f^+$ and $f^-$ are indeed restrictions of some locally defined smooth Lagrangian fibrations (cf. \S \ref{sec:normal}).

\medskip
Let $\pi_{\R}$ be the projection of $\R \times M$ onto $\R$.
Given a point $m \in M$ we study the geometry of a stitched Lagrangian 
fibration $f$ in a neighbourhood of the fibre over $(0,m)$.
For this purpose it is convenient to allow a more general set of coordinates 
on $\R \times M$ than just the smooth ones.

\begin{defi}\label{defi:admissible}  
Let $B$ be a neighbourhood of $(0,m) \in \R \times M$, let  $B^{+} = B \cap  (\R_{\geq 0} \times M)$ and $B^{-} = B \cap  (\R_{\leq 0} \times M)$.  A continuous coordinate chart $(B, \phi)$ around $(0,m)$ is said to be \textbf{admissible} if the components of $\phi=(\phi_1,\ldots ,\phi_n)$ 
satisfy the following properties:
\begin{itemize}
\item[\textit{(i)}] $\phi_1 = \pi_{\R}$;
\item[\textit{(ii)}]  for $j = 2, \ldots, n$ the restrictions of $\phi_j$ to $B^+$ and $B^-$ are locally restrictions of smooth functions on $B$.
\end{itemize}

\end{defi}

\begin{lem} \label{discrep}
Let $f: X \rightarrow \R \times M$ be a stitched Lagrangian fibration and 
let $(B, \phi)$ be an admissible coordinate chart around 
$(0,m) \in \R \times M$. For $j = 2, \ldots , n$, the function
$G_{j}^{\pm} = (\phi_j  \circ  f)|_{f^{-1}(B^{\pm})}$ is the restriction of a $C^\infty$ function on $X$ to $X^\pm$.
Let $\eta_1$ and $\eta_{j}^{\pm}$ be the Hamiltonian vector fields of $\mu$ and $G^{\pm}_{j}$ respectively. Then there are $S^1$ invariant functions $a_j$, 
$j=2, \ldots, n$ on $Z \cap f^{-1}(B)$ such that
\begin{equation} 
      (\eta^{+}_{j} - \eta^{-}_{j})|_{Z \cap f^{-1}(B)} = a_j \, 
                                         \eta_1|_{Z \cap f^{-1}(B)}. 
\label{discrep eq}
\end{equation}
\end{lem}
\begin{proof} Let $\bar{Z} = Z / S^1$, with projection
$p: Z \rightarrow \bar{Z}$ and let $\omega_r$ 
be the Marsden-Weinstein reduced symplectic form on $\bar{Z}$. 
Given a vector field $v$ on $\bar{Z}$, let $\tilde{v}$ be a lift of $v$ on $Z$.
Then we have
\begin{eqnarray*}
    \omega_r(p_{\ast}(\eta_j^+ - \eta_j^-), v) & = &
                         \omega (\eta_j^+ - \eta_j^-, \tilde{v}) \\
        & = & (dG^{+}_{j} -  dG^{-}_{j})(\tilde{v}) = 0, 
\end{eqnarray*}
where the last equality comes from the fact that, being $G$ continuous, 
$G^{+}_{j}|_{Z \cap f^{-1}(B)} = G^{-}_{j}|_{Z \cap f^{-1}(B)}$. Since 
$\omega_r$ is non-degenerate on $Z$, it follows that
                     $$p_{\ast}(\eta_j^+ - \eta_j^-) = 0. $$
Therefore (\ref{discrep eq}) must hold for some function $a_j$, which
 must be $S^1$ invariant since the left-hand side of (\ref{discrep eq})
is $S^1$ invariant. 
\end{proof}

Clearly, when $f$ and the coordinate map $\phi$ are smooth, all the $a_j$'s 
vanish, so equation (\ref{discrep eq}) measures how far $f$ and $\phi$ are 
from being smooth. We will say more about this in the coming sections.

\medskip
We now recall some of the examples which we already introduced and discussed
extensively in \cite{CB-M-torino}.
Consider the following $S^{1}$ action on $\C^3$:
\begin{equation} \label{action}
e^{i\theta}(z_1,z_2,z_3)=(e^{i\theta}z_1,e^{-i\theta}z_2,z_3).
\end{equation}
This action is Hamiltonian with respect to the standard symplectic form $\omega_{\C^3}$. Clearly, it is singular along the surface $\Sigma = \{ z_1=z_2=0 \}$. 
The corresponding moment map is:
\begin{equation}\label{eq mu}
\mu (z_1,z_2,z_3)=\frac{|z_1|^2-|z_2|^2}{2}.
\end{equation}
The only critical value of $\mu$ is $t=0$ and 
$\Crit(\mu )= \Sigma \subset \mu^{-1}(0)$. 

Let $\gamma: \C^2 \rightarrow \C$ be the following piecewise smooth map
\begin{equation}\label{eq. g}
\gamma (z_1,z_2)= \begin{cases}
\frac{z_1z_2}{|z_1|},\quad\text{when}\ \mu \geq 0\\
\\
\frac{z_1z_2}{|z_2|},\quad\text{when}\ \mu <0.
\end{cases}
\end{equation}

In two dimensions we have the following:
\begin{ex} [\textbf{Stitched focus-focus}] \label{broken focus focus}
Consider the map
\begin{equation} \label{fib ff}
f(z_1,z_2)=\left( \frac{|z_1|^2-|z_2|^2}{2}, 
                        \, \log |\gamma (z_1,z_2)+1| \right).
\end{equation}
It is clearly well defined on 
$X = \{ (z_1, z_2) \in \C^2 \ | \ \gamma(z_1, z_2) + 1 \neq 0 \}$ and
it has Lagrangian fibres. 
We showed in \cite{CB-M-torino} that $f$ has the same topology 
of a smooth focus-focus fibration.  The only singular fibre, $f^{-1}(0)$, is a (once) pinched torus. One can easily see that, when restricted to $X - f^{-1}(0)$, $f$ is a stitched 
Lagrangian fibration. The seam is $Z = \mu^{-1}(0) - f^{-1}(0)$. 
Notice that $Z$ has two connected components. Let $\eta_1$ and $\eta_2^{\pm}$ be the 
Hamiltonian vector fields defined as in Lemma~\ref{discrep}. After some computation one 
can verify that
\[ (\eta_2^+ - \eta_2^-)|_{Z} = a \, 
                       \eta_{1}|_Z, \]
where
\[ a = \re  
    \left( \frac{z_1 z_2}{ |z_1|^2 z_1z_2 - |z_1|^3} \right) 
      \left|_{Z} \right.. \]

\end{ex}

There is an analogous model in three dimensions:

\begin{ex} \label{leg}
Consider the map
\begin{equation}\label{eq leg fibr} 
f(z_1, z_2, z_3) = \left( \mu, \, \log |z_3|, \, 
             \log |\gamma (z_1,z_2)-1| \right).
\end{equation}
The discriminant locus of $f$ is 
$\Delta=\{ 0 \} \times \R \times \{0 \} \subset \R^3$. Again, $f$ restricted to $X-f^{-1}(\Delta )$ defines a stitched Lagrangian fibration. 
\end{ex}

\begin{ex} \label{ex amoebous fibr}
Consider the map
\begin{equation} \label{eq. the fibration}
f(z_1, z_2, z_3) = (\mu, \log \frac{1}{\sqrt{2}}| \gamma - z_3|, 
                  \log \frac{1}{\sqrt{2}} |\gamma + z_3 - \sqrt{2}|).
\end{equation}
Let $X$ be the dense open subset of $\C^3$ where $f$ is well defined.
The general construction discussed in \S 5 of \cite{CB-M-torino}
shows that $f$ is a piecewise smooth Lagrangian fibration. 
It contains singular fibres, in fact 
the discriminant locus $\Delta$ of $f$ is depicted in Figure 
\ref{fig: amoeba}. 
\begin{figure}[!ht]
	\centering
	\includegraphics[width=4cm,height=4cm,bb=0 0 258 258]{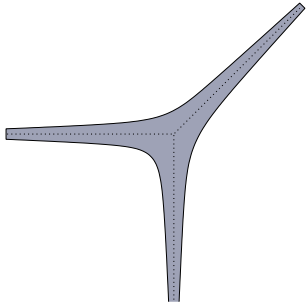}
	\caption{Amoeba of $v_1+v_2+1=0$}
	\label{fig: amoeba}
\end{figure}
One easily checks that $f$ restricted to $X-f^{-1}(\Delta)$
is a stitched Lagrangian fibration. The seam is 
$Z = \mu^{-1}(0) - f^{-1}(\Delta)$, notice that $Z$ has three
connected components. 
Let $\eta_1$ and $\eta_{j}^{\pm}$  be the Hamiltonian vector fields defined as
in Lemma~\ref{discrep}. 
A computation shows that, for $j=2,3$
\[ (\eta_j^{+} - \eta_j^{-})|_{Z} = a_j \, \eta_1 |_{Z}, \]
where
\[ a_2 = - \frac{ \re \left( (\gamma - z_3) 
                  \frac{\zbar_1 \zbar_2}{|z_1|^3} \right)}{|\gamma -z_3|^2 } \]
and
\[ a_3 = - \frac{ \re \left( (\gamma + z_3 - \sqrt{2}) 
                  \frac{\zbar_1 \zbar_2}{|z_1|^3} \right)}{|\gamma + z_3 - \sqrt{2}|^2 }. \]
In \cite{CB-M-torino} we describe the topology of the singular
fibres of $f$ and discuss the relevance of this fibration in the
context of Gross' topological mirror symmetry construction. 
We also show how this example can be perturbed to obtain other interesting
stitched Lagrangian fibrations with discriminant locus of mixed codimension one
and two.
\end{ex}

\section{The first order invariant}
Our goal in this paper is to give a semi-global classification of stitched
Lagrangian fibrations up to \emph{smooth} fibre-preserving symplectomorphism.
For this purpose in this section we restrict our attention to stitched 
Lagrangian fibrations $f: X \rightarrow \R \times M$, where $M = \R^{n-1}$. 
We assume that $f(X)\subseteq\R^n$ is a contractible open neighbourhood $U$ of 
$0 \in \R^n$ and we denote coordinates on $U$ by $b = (b_1, \ldots, b_n)$.
Clearly $(X, f)$ is a topologically trivial torus bundle over $U$. 
Let $U^+ := U \cap \{ b_1 \geq 0 \}$, $U^- := U \cap \{ b_1 \leq 0 \}$
and $\Gamma := U \cap \{ b_1 = 0 \} $.  We assume for simplicity that the 
pair $(U, \Gamma)$ is homeomorphic to the pair $(D^n, D^{n-1})$, 
where $D^n \subset \R^n$ is an $n$-dimensional ball centred at 
$0$ and $D^{n-1} \subset D^n$ is the intersection of $D^n$ with an 
$n-1$ dimensional subspace. 
We have that $X^{\pm} = f^{-1}(U^{\pm})$ and $Z = f^{-1}(\Gamma)$.
If $f_1, \ldots, f_n$ are the components of $f$, then $f_1 = \mu$ is the 
moment map of the $S^1$ action. When $j=2, \ldots, n$,
we denote  $f_j^{\pm} = f_j|_{X^{\pm}}$.  As in Lemma~\ref{discrep} we 
let $\eta_1$ and $\eta_{j}^{\pm}$ denote the Hamiltonian vector fields of
$f_1$ and $f_j^{\pm}$ respectively. Let us recall the notation used in the
proof of Lemma~\ref{discrep}. Since $Z = \mu^{-1}(0)$, the
$S^1$ action on $X$ induces an $S^{1}$-bundle $p: Z \rightarrow \bar{Z}$, where
$\bar{Z} = Z / S^1$. There exists
$\bar{f}: \bar{Z} \rightarrow \Gamma$ such that the following diagram commutes
\[
\xymatrix{
Z \ar[dr]_f \ar[rr]^{p} & & \bar{Z} \ar[dl]^{\bar{f}} \\
&  \Gamma
}
\]
In fact $\bar{f}= (f_2, \ldots, f_n)$, where each component of $\bar f$ is thought of as a function on $\bar{Z}$.
For $b \in \Gamma$ denote by $F_b$ the fibre over $b$ and $\bar{F}_b = F_b\slash S^1$,
clearly $\bar{F}_b = \bar{f}^{-1}(b)$. Denote by
\[ \mathfrak{L} = \ker \bar{f}_{\ast}, \]
the bundle over $\bar Z$ whose fibre at a point $y \in \bar{F}_b$ is $T_y\bar F_b$. 
From Lemma~\ref{discrep} it follows that 
$p_{\ast} \eta_j^{+} = p_{\ast} \eta_j^{-}$, so we can define $\bar{\eta} = (\bar{\eta}_2, \ldots, \bar{\eta}_n)$ to be the frame of $\mathfrak{L}$
where $\bar{\eta}_j = p_{\ast} \eta_j^{\pm}$. 
We say that a section of $\Lambda^k \mathfrak{L}^{\ast}$ is fibrewise 
closed (exact) if it is closed (exact) when viewed as a $k$-form
on each fibre $\bar{F}_b$.  We have the following:

\begin{prop}\label{prop. ell_1}
Let $(X,f)$ be a stitched Lagrangian fibration. If $\ell_{1}$ is the section of $\mathfrak{L}^{\ast}$ defined by
\[ \ell_{1}(\bar{\eta}_j) = a_j, \]
where $a_j$ is the $S^1$-invariant function appearing in (\ref{discrep eq}), 
then $\ell_1$ is fibrewise closed.
\end{prop}
\begin{proof} 
Since $f$ is a Lagrangian submersion, the Hamiltonian vector fields 
$\eta_1, \eta_2^{\pm}, \ldots, \eta_2^{\pm}$ commute and are linearly 
independent. Therefore, for every fixed $b \in \Gamma$, the vector fields 
$\eta_j^{\pm}|_{F_b}$ span $(n-1)$-dimensional integrable distributions $H_{b}^{\pm}$, 
which are horizontal with respect to the $S^1$-bundle $p_b: F_b \rightarrow \bar F_b$.
From the $S^1$ invariance of 
$f_2, \ldots, f_n$, it also follows that $H_{b}^{\pm}$ are $S^1$ invariant.
Thus they define flat connections $\theta_{b}^{\pm}$ of the bundle 
$p_b: F_b \rightarrow \bar F_b$. 
From the properties of flat connections, it follows that 
$\theta_{b}^{-} - \theta_{b}^{+}$ is the pull back of a closed one 
form on $\bar{F}_b$. From 
(\ref{discrep eq}) we obtain
\[ (\theta_{b}^{-} - \theta_{b}^{+})(\eta_{j}^{\pm}) = a_j|_{\bar{F}_b}, \]
i.e. that
\[ \theta_{b}^{-} - \theta_{b}^{+} = p_{b}^{\ast} (\ell_1|_{\bar{F_b}}). \]
Therefore $\ell_1$ is fibrewise closed.
\end{proof}

Clearly, the definition of $\ell_1$ depends on a choice of coordinates on 
$U$. Let $\ell_1^{\prime}$ be a fibrewise closed section of $\mathfrak{L}^{\ast}$.
We say that $\ell_1^{\prime}$ is equivalent to $\ell_1$ up to a change of coordinates on the base if there exists a neighbourhood $W \subseteq U$ of $\Gamma$ and an admissible 
coordinate map $\phi: W \rightarrow \R^n$ such that $\ell_1^{\prime}$ is the 
section associated to $(f^{-1}(W), \phi \circ f)$ via Proposition \ref{prop. ell_1}. Denote by $[\ell_1]$ the class represented by $\ell_1$ modulo this equivalence relation.
We say that a section $\delta$ of $\mathfrak{L}^{\ast}$ is fibrewise 
constant if 
\[ \mathcal{L}_{\bar{\eta}_j} \, \delta = 0, \]
for all $j=2, \ldots, n$, here $\mathcal{L}_{\bar{\eta}_j}$ denotes the Lie derivative. 
One can easily check that the latter definition is 
independent of the admissible coordinates on the base used to define
$\bar{\eta}_j$. We have the following

\begin{prop} \label{base:chcoor}
A section $\ell_1^{\prime}$ of $\mathfrak{L}^{\ast}$ is equivalent to 
$\ell_1$ up to a change of coordinates on the base if and only if 
\[
\ell_1^{\prime} = \ell_1 + \delta,
\]
where $\delta$ is fibrewise constant. In particular, the class of $\ell_1$ may be written as
\[
[\ell_1]=\{ \ell_1+\delta\mid\delta\ \textrm{is fibrewise constant} \}.
\]
\end{prop}
\begin{proof} Given an admissible change of coordinates 
$\phi: W \rightarrow \R^n$,
we must have $\phi_1 = b_1$. Moreover the partial derivatives 
$\partial_{k} \phi_j$ are defined and continuous on $W$ for all 
$k,j = 2, \ldots, n$. As far as derivatives with respect to $b_1$ are
concerned, only left and right derivatives are defined and 
smooth on $\Gamma$, i.e. only $\partial_1 \phi^+_j$ and $\partial_1 \phi^-_j$,
which may a priori differ. Let $(\eta_j^{\prime})^{\pm}$ be the Hamiltonian
vector fields on $Z$ corresponding to $\phi_j \circ f$, with 
$j =2, \ldots, n$, and let $\bar{\eta}^{\prime}_j = 
p_{\ast}(\eta_j^{\prime})^{\pm}$ . An easy calculation shows that
\[ (\eta_j^{\prime})^{\pm} = \partial_1 \phi^{\pm}_j \, \eta_1 + 
                     \sum_{k=2}^{n} \partial_k \phi_j \, \eta^{\pm}_k. \]
In particular this implies
\begin{equation} \label{chng} 
       \bar{\eta}^{\prime}_j  = \sum_{k=2}^{n} \partial_k \phi_j \, \bar{\eta}_k. 
\end{equation}
and
\begin{eqnarray*}
 (\eta_j^{\prime})^{+}- (\eta_j^{\prime})^{-} & = &
 (\partial_1 \phi^{+}_j -  \partial_1 \phi^{-}_j) \, \eta_1 + 
             \sum_{k=2}^{n} \partial_k \phi_j \, (\eta^{+}_k - \eta^{-}_k) \\
          & = & (\partial_1 \phi^{+}_j -  \partial_1 \phi^{-}_j + 
             \sum_{k=2}^{n} a_k\partial_k \phi_j ) \, \eta_1. 
\end{eqnarray*}
If $\ell_1^{\prime}$ is the 1-form associated to $\phi \circ f$ via Proposition \ref{prop. ell_1}, then
by definition we must have
\[ \ell_1^{\prime}(\bar{\eta}^{\prime}_j) = 
  \partial_1 \phi^{+}_j -  \partial_1 \phi^{-}_j + 
             \sum_{k=2}^{n} a_k\partial_k \phi_j. \]
Let $\delta$ be the section of $\mathfrak{L}^{\ast}$ defined by
\begin{equation} \label{base:discr}
  \delta(\bar{\eta}^{\prime}_j) =  
                 \partial_1 \phi^{+}_j -  \partial_1 \phi^{-}_j. 
\end{equation}
Then, also using (\ref{chng}), we see that
\[ (\ell_1 + \delta)(\bar{\eta}^{\prime}_j) =  
                   \ell_1^{\prime}(\bar{\eta}^{\prime}_j). 
                                                                             \]
Moreover, from (\ref{base:discr}) one can see that $\delta (\bar\eta'_{j})$ descends to a function on $\Gamma$ and therefore $\delta$ is fibrewise constant.

\medskip
Now suppose that $\delta$ is a fibrewise constant section of $\mathfrak{L}^{\ast}$. 
Let
\[ \delta(\bar{\eta}_j) = d_j. \]
Since $\delta$ is fibrewise constant, the $d_j$'s are fibrewise constant 
functions on $\bar{Z}$, i.e. they descend to functions on $\Gamma$. 
Define the following map
\[\phi(b_1, \ldots, b_n) = \begin{cases}
(b_1, \, b_2 + d_2(b_2, \ldots, b_n) \, b_1, \ldots, \,  b_n + 
          d_n(b_2, \ldots, b_n) \, b_1),\quad\text{when}\ b_1 \geq 0\\
\\
\I,\quad\text{when}\ b_1 <0.
\end{cases} \]
It is a well defined admissible coordinate map on some open neighbourhood
of $\Gamma$. It is also clear that (\ref{base:discr}) holds. 
\end{proof}

\begin{defi} We call $\ell_1$ the \textbf{first order invariant} of the stitched fibration
$(X, \omega ,f)$. 
\end{defi}
The name ``invariant" in the above Definition will be fully justified later on. 

\medskip
It is clear from the proof of Proposition \ref{base:chcoor} that $\delta$ is a first order measure of how far the change of coordinates on the base is from being smooth, in particular if it is smooth then $\delta = 0$. 

We also have the following:
\begin{cor}
If there exists an admissible change of coordinates on the base
which makes the stitched Lagrangian fibration smooth, then $\ell_1$ is 
fibrewise constant.
\end{cor} 
\begin{proof} It is clear that if $\phi \circ f$ is smooth then we 
must have that its first order invariant $\ell_1^{\prime}$ is zero.
It then follows from Proposition~\ref{base:chcoor} that $\ell_1$ must
be fibrewise constant.
\end{proof}

We now describe action-angle coordinates of a stitched Lagrangian fibration $f: X \rightarrow U$. Let $\alpha$ be a 1-form on $U$. Since $f^\pm$ is the restriction of a smooth map, $\alpha$ pulls back to an honest smooth 1-form $\alpha^\pm$ defined on a neighbourhood of $Z$. The latter defines a smooth vector field $v_\alpha^\pm$ determined by the equation (\ref{eq contr.}). The flow of $v_\alpha^\pm$, when restricted to $X^\pm$, is fibre-preserving. This induces an action of $T^\ast_bU$ on the fibre $(f^\pm)^{-1}(b)$ for all $b\in U^\pm$. Let $\sigma: U \rightarrow X$ be a continuous section which is smooth and Lagrangian when restricted to $U^{\pm}$. Then, as explained in \S \ref{section aa}, there is a maximal smooth lattice $\Lambda_{\pm} \subset T^{\ast}U^{\pm}$ and a diagram
\[
\xymatrix{
T^{\ast}U^{\pm} / \Lambda_{\pm}  \ar[dr]_{\pi^\pm} \ar[rr]^{\Theta^\pm} & & X^\pm \ar[dl]^{f^\pm} \\
&  U^\pm
}
\]
where $\Theta^{\pm}$ is a symplectomorphism and $\pi^{\pm}$ is
the standard projection. Let $\Phi_{\eta_1}^{t}, \Phi_{\eta_2^{\pm}}^{t} \ldots, \Phi_{\eta_n^{\pm}}^{t}$ denote the flow of $\eta_1, \eta_2^{\pm}, \ldots, \eta_n^{\pm}$ 
respectively. Then
\begin{equation} \label{theta:flows}
               \Theta^{\pm}: (b, \, \sum_{j} t_j \, db_j) \mapsto  
          \Phi_{\eta_1}^{t_1} \circ \Phi_{\eta_2^{\pm}}^{t_2} \circ 
                        \ldots \circ \Phi_{\eta_n^{\pm}}^{t_n}(\sigma(b)),
\end{equation}
and
\[ \Lambda_{\pm} = \{ (b, \, \sum_{j} T_j \, db_j) \in T^{\ast}U^{\pm} \ | \ 
 \Phi_{\eta_1}^{T_1} \circ \Phi_{\eta_2^{\pm}}^{T_2} \circ \ldots 
                      \circ \Phi_{\eta_n^{\pm}}^{T_n}(\sigma(b))
          = \sigma(b) \} \]
Now let 
\[ \lambda_1= db_1. \]
The $S^1$ action implies $db_1 \in \Lambda_{\pm}$.
Let us denote a basis for $\Lambda_{\pm}$ by 
$\{\lambda_1, \lambda_2^{\pm}, \ldots, \lambda_n^{\pm} \}$, where
\[ \lambda_{j}^{\pm} = \sum_{k=1}^{n} T_{jk}^{\pm} db_k. \]
The $S^1$ action on $X$ corresponds to translations along
the $\lambda_1$ direction. Let
\[ Z^{\pm} = (\pi^{\pm})^{-1}(\Gamma). \]
If we denote
\[ \bar{\lambda}_{j}^{\pm} = \lambda_{j}^{\pm} \mod db_1 \]
and let $\bar{\Lambda}^{\pm} = 
       \spn \inner{ \bar{\lambda}_{2}^{\pm}}{, \ldots, \bar{\lambda}_{n}^{\pm}}_{\Z}$,
then $\Theta^{\pm}$ identifies $Z / S^1$ with 
\[ \bar{Z}^{\pm} = T^{\ast} \Gamma / \bar{\Lambda}_{\pm}. \]
Denote by $\bar{t} = (t_2, \ldots, t_n)$ the coordinates on 
the fibres of $\bar{Z}^-$.

\medskip
Now observe that, due to the discrepancy (\ref{discrep eq}) between $\eta_{j}^{+}$
and $\eta_{j}^{-}$ along $Z$,  $\Theta^+$ and $\Theta^-$ behave 
differently on fibres lying over $\Gamma$. We have the diagram:
\[
\xymatrix{
Z^- \ar[dr]_{\Theta^-}  & & Z^+ \ar[dl]^{\Theta^+} \\
&  Z
}
\]
and the difference between the two maps is measured by
\[ (\Theta^{+})^{-1} \circ \Theta^{-} :   
       Z^{-} \rightarrow  Z^{+}. \]
We have the following characterisation of this map:
\begin{prop} \label{map:Q}
The discrepancy (\ref{discrep eq}) between the Hamiltonian vector fields
of the stitched Lagrangian fibration $f: X \rightarrow U$ induces the map
\[ Q = (\Theta^{+})^{-1} \circ \Theta^{-} \]
between $Z^-$ and $Z^+$. Let
\[ \ell_1^- = (\Theta^-)^{\ast} \ell_1. \]
Then, computed explicitly in the canonical coordinates
on  $ T^{\ast} U^{-}$ and $T^{\ast} U^{+}$, $Q$ is given by
\begin{equation} \label{glue:Q}
 Q: 
(b, t_1, \bar{t}) \mapsto \left(b, \  t_1 - 
    \int_0^{\bar{t}} \ell_1^-, \ \bar{t} \right), 
\end{equation}
where $(b, \bar{t})$ are the canonical coordinates on $\bar{Z}^-$ and 
the integral is a line integral in $T^{\ast}_b\Gamma$
along a path joining $(b,0)$ and $(b, \bar t)$.
\end{prop}
\begin{proof}  Let $(b_1, \ldots, b_n, t_1, \ldots, t_n)$ and 
$(b_1, \ldots, b_n, y_1, \ldots, y_n)$ be the canonical coordinates on  
$ T^{\ast} U^{-}$ and $T^{\ast} U^{+}$ respectively.
From its definition, we see that $\Theta^{+}$ 
identifies $\eta_1, \eta_{2}^{+}, \ldots, \eta_{n}^{+}$  with 
$\partial_{y_1}, \ldots, \partial_{y_n}$ and w.l.o.g. we can assume that it sends $\sigma$ to the zero section of $T^\ast U$. Therefore (\ref{discrep eq}) becomes
\begin{equation} \label{vf:coo}
  \eta_{j}^{-} = \partial_{y_j} - (a_j \circ \Theta^{+})
                                                \ \partial_{y_1}.
\end{equation}
Notice that $a_j \circ \Theta^{+}$ is independent of $y_1$.
Computing the flows of $\eta_1, \eta_2^{-}, \ldots, \eta_n^{-}$ in these coordinates is
not difficult and it turns out that $Q$ is given by
\[ Q:
(b, t_1, \ldots, t_n) \mapsto \left( b, \ t_1 - \sum_{j=2}^{n} 
    \int_0^{t_j} a_j \circ \Theta^- (b, t_2, \ldots, t_{j-1}, t,0, \ldots, 0) dt, 
\ t_2, \ldots, \ t_n \right), \]
which is equivalent to (\ref{glue:Q}), since $\ell_1$ is fibrewise closed\footnote{To verify that the above expression of $Q$ is correct, it is enough to check that $\partial_{t_j}Q=\eta^-_j\circ Q$.}.
\end{proof}

We now explain how the map $Q$ matches the periods in $\Lambda^-$
with those in $\Lambda^+$. The maps $\Theta^{\pm}$ naturally identify
$\Lambda_{\pm}$ with $H_1( X, \Z) \cong \Z^n$, but in general $\Theta^{-}$ does
it differently from $\Theta^+$. Let $\gamma_1$ be the cycle represented by the
orbit of the $S^1$ action. We know that $\gamma_1$ always corresponds to
the period $db_1$. 

\medskip
We have the following
\begin{cor} \label{broken:per} Suppose we choose bases  
$\{ \lambda_{1}, \lambda_2^{\pm}, \ldots, \lambda_n^{\pm} \}$ of 
$\Lambda_{\pm}$ corresponding to two bases
$\gamma^{\pm} = \{ \gamma_1,\gamma_2^{\pm}, \ldots, \gamma_n^{\pm} \}$ of $H_1( X, \Z)$, 
such that
\begin{itemize}
\item[\textit{(i)}] $\gamma_1$ is represented by an orbit of the $S^1$ action,
\item[\textit{(ii)}] $\gamma_j^{+} = \gamma_j^{-} + m_j \gamma_1$, for some $m_2, \ldots,m_n\in\Z$.
\end{itemize}
Then at a point $b \in \Gamma$ we have 
  \begin{equation} \label{per:corr}
    \lambda_{j}^{+}(b) = \lambda_{j}^{-}(b) + \left( m_j -
         \int_{\bar{\lambda}_{j}^{-}} \ell_1^{-}  
                                               \right) \, \lambda_1, 
  \end{equation}
where the integral of $\ell_1^-$ is taken along the cycle represented
by $\bar{\lambda}_{j}^{-}$.
In particular
\begin{equation} \label{p2:id}
  \bar{\lambda}_{j}^{+}(b) = \bar{\lambda}_{j}^{-} (b).
\end{equation}
\end{cor}
\begin{proof}
To obtain (\ref{per:corr}) it suffices to observe that since
$m_j \lambda_1 + \lambda_{j}^{-}$ and $\lambda_{j}^{+}$ have to represent the same
$1$-cycle in $f^{-1}(b)$, they must be mapped one to
the other by  $Q$. The
result is therefore obtained by applying (\ref{glue:Q}) to $m_j\lambda_1+\lambda_j^-$.
\end{proof}
\begin{rem}
Condition $(ii)$ means that under the map $p_*:  H_1(X,\Z ) \rightarrow 
H_1(X\slash S^1,\Z )$, bases $\gamma^+$ and $\gamma^-$ are mapped to 
the same base of $H_1(X\slash S^1,\Z )$.
We will need to consider condition $\textit{(ii)}$ in 
 \S~\ref{sec. stitch w/mono}  where we discuss stitched Lagrangian fibrations over non simply connected bases, for which non-trivial monodromy may occur. 
\end{rem}
\begin{rem} \label{quotient}
From Proposition~\ref{map:Q} and Corollary \ref{broken:per} it follows that 
$\bar{\Lambda}_- = \bar{\Lambda}_+$ and that on the quotients $\bar{Z}^+$
and $\bar{Z}^-$, $Q$ acts as the identity. Therefore, if we let
\[ \ell_1^+ = (\Theta^+)^{\ast} \ell_1, \]
then we have $\bar{Z}^+ = \bar{Z}^-$ and $\ell_1^+ = \ell_1^-$. Thus we
can remove the $+$ and $-$ signs and denote
\[ \bar{\lambda}_j =  \bar{\lambda}_{j}^{+} = \bar{\lambda}_{j}^{-} \]
\[ \bar{\Lambda} = \bar{\Lambda}_- = \bar{\Lambda}_+ \]
and, with slight abuse of notation, identify $\bar{Z}$  with 
$T^{\ast} \Gamma / \bar{\Lambda}$ and $\bar{f}$ with
the projection $\bar{\pi}: \bar{Z} \rightarrow \Gamma$. Notice then that
$\mathfrak{L}$ is identified with $\ker \bar{\pi}_{\ast}$ and $\ell_1$
with $\ell_1^{\pm}$. 
\end{rem}

It is natural to consider bases of $H_1(X, \Z)$ satisfying conditions
$(i)$ and $(ii)$ also because of the following
\begin{lem} \label{broken:action}  
Let 
$\{ \gamma_1,\gamma_2^{\pm}, \ldots, \gamma_n^{\pm} \}$ be bases of 
$H_1( X, \Z)$
satisfying conditions $(i)$ and $(ii)$ of Corollary~\ref{broken:per}
and let $\alpha^{\pm}: U^{\pm} \rightarrow \R^n$ be the corresponding 
action coordinates satisfying $\alpha^{\pm}(0) = 0$. Then the map
\begin{equation}\label{eq. stitched action}
\alpha = \begin{cases}  \alpha^{+} \quad\text{on} \ U^{+}, \\

                           \alpha^{-} \quad\text{on} \ U^{-}, 
                \end{cases}
\end{equation}
is an admissible change of coordinates.
\end{lem}
\begin{proof} 
Action coordinates 
$\alpha^{\pm} = (\alpha_1^\pm ,\ldots , \alpha_n^\pm )$ are defined 
by the integral
\[
\alpha^\pm_j(b) = \int_{0}^{b} \lambda_j^{\pm}.
\]
along a curve in $U^\pm$ joining $0$ and $b$. When $j=1$, this gives
 $\alpha_1^+=\alpha_1^-=b_1$.
Clearly $\alpha$ is a diffeomorphism when restricted to $U^+$ or $U^{-}$.
Moreover $\alpha$ is injective. 
The fact that $\alpha^+$ and $\alpha^-$ coincide along $\Gamma$ follows from 
(\ref{p2:id}) and the connectedness of $\Gamma$. In fact (\ref{p2:id})
 implies that when $b \in \Gamma$ the above integral gives
\[
\alpha^+_j(b) = \int_{0}^{b} \bar{\lambda}_j^{+}= 
            \int_{0}^{b} \bar{\lambda}_j^{-} = \alpha^-_j(b).
\]
This concludes the proof.
\end{proof} 

\begin{rem} \label{rem:action}
The upshot of Lemma \ref{broken:action} is that after a change of coordinates as in (\ref{eq. stitched action}) we can always assume that the coordinates on
the base $U$, when restricted to $U^{\pm}$, are action coordinates corresponding
to bases $\{ \gamma_1,\gamma_2^{\pm}, \ldots, \gamma_n^{\pm} \}$ of $H_1( X, \Z)$
satisfying $(i)$ and $(ii)$ of Corollary~\ref{broken:per}. Then 
$\{ db_1, db_2, \ldots, db_n \}$ form a basis of $\Lambda_{+}$ and $\Lambda_-$. From 
(\ref{per:corr}) it also follows that, in view of the identifications of
Remark~\ref{quotient}, $\ell_1$ must satisfy
\[  \int_{\bar{\lambda}_{j}} \ell_1 = m_j. \]
The reader should be warned at this point that, although the map $\alpha$ as in (\ref{eq. stitched action}) allows us to find action coordinates on both $U^+$ and $U^-$,  we still have two different sets of action-angle coordinates, $(b_1,\ldots ,b_n,y_1,\ldots ,y_n)$ on $X^+$ and $(b_1,\ldots ,b_n,t_1,\ldots ,t_n)$ on $X^-$. This is due to the discrepancy between $\Theta^+$ and $\Theta^-$, which makes the map:
\[
\Theta = \begin{cases}  (\Theta^{+})^{-1} \quad\text{on} \ X^{+} \\

                           (\Theta^{-})^{-1} \quad\text{on} \ X^{-} 
                \end{cases}
\]
discontinuous along the seam $Z$. As pointed out before, this discrepancy is measured by $\ell_1$.
\end{rem}

In the next theorem we show that any fibrewise closed section 
$\ell_1 \in \mathfrak{L}^{\ast}$ can be the first 
order invariant of a stitched Lagrangian fibration.

\begin{thm} \label{broken:constr}
Let $U$ be an open contractible neighbourhood of $0 \in \R^n$ such that
$\Gamma = U \cap \{ b_1 = 0 \}$ is contractible. Let
$\bar{\Lambda} \subseteq T^{\ast} \Gamma$ be the lattice spanned by
$\{ db_2, \ldots, db_n \}$, and let
$\bar{Z} = T^{\ast} \Gamma / \bar{\Lambda}$, with projection
$\bar{\pi}: \bar{Z} \rightarrow \Gamma$ and bundle
$\mathfrak{L} = \ker \bar{\pi}_{\ast}$. Given integers $m_2, \ldots, m_n$
and a
smooth, fibrewise closed section $\ell_{1}$ of $\mathfrak{L}^{\ast}$
such that
\begin{equation} \label{int:cond}
  \int_{db_j} \ell_1 = m_j \ \ \ \ \text{for all} \ \ j=2, \ldots, n,
\end{equation}
there exists a smooth symplectic manifold $(X, \omega)$ and a stitched
Lagrangian fibration $f: X \rightarrow U$ satisfying the following
properties:
\begin{itemize}
\item[\textit{(i)}] the coordinates $(b_1, \ldots, b_n)$ on $U$
          are action coordinates of $f$ with $\mu = f^{\ast}b_1$;
    \item[\textit{(ii)}] the periods $\{ db_1, \ldots, db_n \}$,
restricted to $U^{\pm}$
          correspond to basis $\{ \gamma_1, \gamma_2^{\pm}, \ldots,
          \gamma_n^{\pm} \}$ of $H_1(X, \Z)$ satisfying $(i)$ and $(ii)$ of
          Corollary~\ref{broken:per};
    \item[\textit{(iii)}]  $\ell_1$ is the first order invariant of $(X,f)$.
\end{itemize}

\end{thm}
\begin{proof} We regard the two halves of $U$, $U^+$ and
$U^-$ defined as before, as disjoint sets. Let $\Lambda_{\pm}$ be the lattices in
$T^{\ast}U^{\pm}$
spanned by $\{ db_1, db_2, \ldots, db_n \}$ and define
$X^{\pm} = T^{\ast}U^{\pm} / \Lambda_{\pm}$, with corresponding
projections $\pi^{\pm}$.
Let $Z^{\pm} = \partial X^{\pm} = (\pi^{\pm})^{-1}(\Gamma)$.
Translations along the $db_1$ direction define an $S^1$ action on
$Z^{\pm}$ such
that $\bar{Z} = Z^{\pm} / S^1$.
On $T^{\ast}U^{+}$ and $T^{\ast}U^{-}$, we consider canonical
coordinates $(b_1, \ldots, b_n, y_1, \ldots, y_n)$ and
$(b_1, \ldots, b_n, t_1, \ldots, t_n)$ respectively (or $(b,y)$ and
$(b,t)$ for short). Coordinates on $\bar{Z}$ are given by $(b, \bar{t})$
(or $(b, \bar{y})$), where  $b=(0, b_2, \ldots, b_n) \in \Gamma$ and
$\bar{t} = (t_2, \ldots, t_n)$ (or $\bar{y} = (y_2, \ldots, y_n)$).
For $j =2, \ldots, n$, let
\[ a_{j} = \ell_1(\partial_{t_j})=  \ell_1(\partial_{y_j})\]
On $X^+$, let $\eta_1 = \partial_{y_1}$ and
$\eta_j^{+} = \partial_{y_j}$ then on $Z^{+}$ we can define
vector fields
\[ \eta_{j}^{-} = \eta_j^{+} - a_j \ \eta_1, \]
which is coherent with (\ref{vf:coo}).  We can define a map
$Q: Z^- \rightarrow Z^+$ by composition of the
flows of $\eta_1, \eta_2^{-}, \ldots, \eta_n^{-}$, i.e.
\[ Q: (b, t_1, \ldots, t_n) \mapsto \Phi_{\eta_1}^{t_1}
                    \circ \Phi_{\eta_{2}^{-}}^{t_2} \circ \ldots \circ
       \Phi_{\eta_{n}^{-}}^{t_n} (b,0). \]
Clearly, $Q$ can be written as in (\ref{glue:Q}).
One can easily see that the properties of $\ell_1$ ensure that $Q$ is
a well defined fibre-preserving diffeomorphism which sends the cycles
represented by $db_1$ and $db_j$ in $H_1(Z^-,\Z )$ to the cycles
represented by $db_1$ and $db_j-m_jdb_1$ in $H_1(Z^+,\Z )$, $j = 2,
\ldots, n$, respectively. Intuitively, $Q$ identifies fibres of $\pi^-$
inside $Z^-$ with fibres of $\pi^+$ inside $Z^+$ after the latter ones
have been twisted by iteratively flowing in the direction of $\eta_j^-$,
$j=2,\ldots ,n$. Topologically we define
\[ X = X^+ \cup_{Q} X^- .\]
To give $X$ smooth and symplectic structures we have to extend
the gluing map $Q$ to open neighbourhoods of $Z^+$ and $Z^-$.
Let open sets $\tilde{U}^+$ and $\tilde{U}^{-}$ be small
enlargements of $U^+$ and $U^-$ respectively, obtained by
joining small open neighbourhoods of $\Gamma$ to $U^+$ and
$U^-$. Extend $\Lambda_{\pm}$ to lattices of
$T^{\ast}\tilde{U}^{\pm}$ in a constant way. We look for neighbourhoods
$V^{\pm}$ of $Z^{\pm}$ inside $T^{\ast}\tilde{U}^{\pm}/ \Lambda_{\pm}$
and a symplectomorphism $\tilde{Q}: V^- \rightarrow V^+$ extending $Q$.
One can achieve this by considering an ``auxiliary" fibration. Suppose for
now that we could find a neighbourhood $V^{+}$
of $Z^+$ and a smooth, proper $S^1$-invariant Lagrangian fibration
$u: V^+ \rightarrow \R^n$, with components $u_j$ such that:
\begin{equation} \label{ext:ham}
\begin{array}{ll}
u_1=b_1, & \\
u|_{Z^+} = \pi^+, & \\
\eta_{u_j}|_{Z^+} = \eta_{j}^{-}, &\textrm{when}\ j=2, \ldots, n.
\end{array}
\end{equation}
This amounts to prescribing zero and first order terms of $u$ along $Z^+$ in 
the Taylor expansion of $u$ with respect to $b_1$.
Now inside $\tilde{U}^-$ there will be a small open neighbourhood $W$ of
$\Gamma$ and a symplectomorphism:
\[ \tilde{Q}: V^- \rightarrow V^+, \]
where $V^-:=(\pi^-)^{-1}(W)$ and

 \[ \tilde Q: (b, t_1, \ldots, t_n) \mapsto \Phi_{\eta_1}^{t_1}
                    \circ \Phi_{\eta_{u_2}}^{t_2} \circ \ldots \circ
       \Phi_{\eta_{u_n}}^{t_n} (b,0). \]

In other words, $\tilde Q$ is the action-angle coordinate map associated
to the fibration $u: V^+ \rightarrow \R^n$, computed with respect to the
cycles $\{ db_1, -m_2 db_1 + db_2, \ldots, -m_n db_1 + db_n \}$ (it may be
necessary, for this purpose, to restrict to a smaller $V^+$). From
(\ref{ext:ham}) it follows that $\tilde Q$ extends $Q$. We
define
\[ X = (X^+ \cup V^+) \cup_{\tilde{Q}} (X^- \cup V^-). \]
and the stitched Lagrangian fibration to be
\[ f =  \begin{cases} \pi^+ \quad\text{on} \ X^+ \\
                      \pi^- \quad\text{on} \ X^-.
        \end{cases}   \]
Due to the non-triviality of the gluing map $\tilde{Q}$ used to
define $X$, $f$ is in general piecewise smooth.
In fact if we pull back $f$ via the inclusion
$ X^+ \cup V^+ \hookrightarrow X$, then we obtain
\[ f|_{ X^+ \cup V^+} =  \begin{cases} \pi^+ \quad\text{on} \ b_1 \geq 0 \\
                      u \quad\text{on} \ b_1 \leq 0.
        \end{cases}   \]
This is because $\pi^- = \tilde{Q}^{\ast} u$.
By construction $(X, \omega)$ and $f$ satisfy the conditions $(i) - (iii)$.

\medskip
Now we prove that a fibration $u:V^+\rightarrow\R^n$ satisfying
(\ref{ext:ham}) exists.
For every $b \in \Gamma$,  consider the following one-parameter family of
closed
$1$-forms on the fibre $F_{b} = (\pi^+)^{-1}(b)$
\[ \ell(r) = r(dy_1 + \ell_1), \]
where $r \in \R$. For every $r$, the graph of $\ell(r)$ defines a Lagrangian
submanifold inside $T^{\ast}F_b$. For $r$ sufficiently small,
let $L_{r,b}$ be the image of the graph of $\ell(r)$ under the
symplectomorphism
\[ (y_1, \ldots, y_n, \sum_{k=1}^{n} x_k dy_k)
            \mapsto (x_1, b_2 + x_2, \ldots, b_n + x_n, y_1, \ldots, y_n),\]
between a neighbourhood of the zero section of $T^{\ast}F_b$ and a
neighbourhood
of $F_b$ inside $T^{\ast}\tilde{U}^{+}/ \Lambda_{+}$. Then there will be a
sufficiently small neighbourhood $V^+$ of $Z^+$ which is fibred by the
submanifolds $L_{r,b}$, i.e. on which the manifolds $L_{r,b}$ are the fibres
of a Lagrangian fibration $u: V^+ \rightarrow \R^n$.
This is due to the fact that the map
\[ (r, b_2, \ldots, b_n, y_1, \ldots, y_n)
              \mapsto (r, b_2 + r a_2(b,\bar{y}), \ldots, b_n + r
a_n(b,\bar{y}),
                                 y_1, \ldots, y_n) \]
is a diffeomorphism when restricted to a neighbourhood of 
$\{0 \} \times Z^+$ inside $\R \times Z^+$. We now show that a possible choice of $u$ also
satisfies (\ref{ext:ham}). Notice that $u$ will be $S^1$-invariant since
its fibres $L_{r,b}$ are $S^1$-invariant.
Given $(b^{\prime}, y^{\prime}) \in V^+$, there exists a unique
$(r,b) \in \R \times Z^+$ such that $L_{r,b} \subset V^+$ and
$(b^{\prime}, y^{\prime}) \in L_{r,b}$.
In fact $(r,b)$ can be determined
as a function of $(b^{\prime}, y^{\prime})$ by solving the non linear system
\begin{equation} \label{sistema}
     \begin{cases} r = b_1^{\prime} \\
                 b_j + r a_j(b, y^{\prime}) = b_j^{\prime} \quad\text{when}
                              \ j=2, \ldots, n
   \end{cases}
\end{equation}
using the implicit function theorem.
Now define
\[ u_1(b^{\prime}, y^{\prime}) = b_1^{\prime} \]
and, when $j=2, \ldots, n$
\begin{equation} \label{ugei}
         u_j(b^{\prime}, y^{\prime}) = b_j,
\end{equation}
where $b_j$ (and thus $b$) are functions of $(b^{\prime}, y^{\prime})$.
Notice
that $S^1$-invariance of $u_j$ can also be seen from the fact that $u_j$ is
independent of $y_1$. It is clear that, when $j = 2, \ldots, n$
\[ \begin{cases}
    \partial_{y_k^{\prime}} u_j |_{Z^+} = 0 \quad\text{for all} \
                   k = 1, \ldots,n \\
   \partial_{b_k^{\prime}} u_j |_{Z^+} = \delta_{kj}
\quad\text{for all}
                             \ k=2, \ldots, n.
     \end{cases} \]
Therefore
      \[ \eta_{u_j}|_{Z^+}     =    \partial_{b_1^{\prime}} b_j \, \partial_{y_1} +
\partial_{y_j}. \]
Using (\ref{sistema}) we compute that
\[ \partial_{b_1^{\prime}} b_j|_{Z^+} = - a_j, \]
which proves that conditions (\ref{ext:ham}) are satisfied.
\end{proof}

\section{Higher order terms} \label{sec:ho-inv}
In Theorem~\ref{broken:constr} we provided a (local) construction of stitched Lagrangian 
fibrations with any given first order invariant satisfying integrality
conditions (\ref{int:cond}). It involved the choice of a Poisson commuting 
set of
functions $u_1, \ldots, u_n$ (producing a Lagrangian fibration $u$) defined
on a neighbourhood of $Z$ and with prescribed $0$-th and $1$-st order terms (cf. (\ref{ext:ham})). 
In general there may be many choices of such functions giving stitched 
Lagrangian fibrations which are not fibrewise symplectomorphic. It is
necessary to look at higher order terms. In this Section we give a description of
these higher order terms and prove an existence result of stitched Lagrangian 
fibrations with prescribed higher order terms.. 

\medskip
We fix here some basic notation. Let $(b_1, \ldots, b_n)$ be standard
coordinates on $\R^n$. Let $\R^{n-1}$ be embedded in $\R^n$ as the subset
$\{ b_1 = 0 \}$ and let $\Gamma \subset \R^{n-1}$ be an open
neighbourhood of $0 \in \R^{n-1}$. We will denote by $U$ an
open neighbourhood of $\Gamma$ in $\R^n$. 
We assume that the pair $(U, \Gamma)$ is 
diffeomorphic to the pair $(D^{n}, D^{n-1})$ where $D^k \subset \R^k$
is a unit ball centred at $0$.
Denote $U^+ = U \cap \{ b_1 \geq 0 \}$ and 
$U^- = U \cap \{ b_1 \leq 0 \}$. Then $\Gamma = U^+ \cap U^-$.
Let $\Lambda$ be the lattice in $T^{\ast} U$ generated
by $\{ db_1, \ldots, db_n \}$ and consider $T^{\ast}U / \Lambda$
with the standard symplectic form and with projection onto $U$ denoted
by $\pi$. 
We assume $S^1$ acts on $T^{\ast}U / \Lambda$ via translations along the $db_1$ direction.
Let $Z = \pi^{-1}(\Gamma)$ and $\bar{Z} = Z / S^1$. If $\bar{\Lambda}$ denotes
the lattice in $T^{\ast} \Gamma$ spanned by $\{ db_2, \ldots, db_n \}$,
we have $\bar{Z} = T^{\ast} \Gamma / \bar{\Lambda}$ with projection $\bar{\pi}$. 
Given $b \in \Gamma$, we denote $F_b = \pi^{-1}(b)$ and 
$\bar{F}_b = \bar{\pi}^{-1}(b) = F_b / S^1$.
Canonical coordinates on $T^{\ast} U$ are denoted by
$(b,y)=(b_1, \ldots, b_n, y_1, \ldots, y_n)$. We also have the bundle
$\mathfrak{L} = \ker \bar{\pi}_{\ast}$. 

Throughout this section we will study the set defined in the following
\begin{defi}
We define $\mathscr U_{\bar{Z}}$ to be the set of pairs $(V,u)$
where $V$ is a neighbourhood of $Z$ and $u: V \rightarrow \R^n$ is a
 $C^\infty$, proper, $S^1$-invariant, Lagrangian submersion, with components
 $(u_1, \ldots, u_n)$, such that $u|_{Z} = \pi$ and $u_1 = b_1$.
\end{defi}

Given $(V,u) \in \mathscr U_{\bar{Z}}$,  let $Y^+ := \pi^{-1}(U^+)$, 
$Y:= Y^+ \cup V$, $Y^- := Y \cap \pi^{-1}(U^-)$ and define the map 
$f_u: Y \rightarrow \R^n$ by
\begin{equation} \label{u:st}
       f_u = \begin{cases}
                u \quad\text{on} \  Y^-, \\
               \pi \quad\text{on} \ Y^+.
       \end{cases} 
\end{equation}
Clearly $(Y,f_u)$ is a stitched Lagrangian fibration. We study the aforementioned higher order terms
of such fibrations. 

\begin{prop} \label{u:tay eq} 
Let $(V,u) \in \mathscr U_{\bar{Z}}$.
For every $N \in \N$ and $j=2, \ldots, n$, 
consider the $N$-th order Taylor series expansion 
of $u_j$ in the variable $b_1$, evaluated at $b_1=0$:
\begin{equation} \label{ugei:tay}
 u_j = \sum_{k=0}^{N} S_{j,k} b_1^k + o(b_1^N), 
\end{equation}
where $S_{j,k}$ are smooth functions on $Z$ which are $S^1$ invariant (i.e. independent of $y_1$). 
For every $m \in \N$, define the following sections of $\mathfrak{L}^{\ast}$ 
and $\Lambda^2 \, \mathfrak{L}^{\ast}$ respectively
\begin{equation}\label{eq:S_m}
S_m = \sum_{j=2}^{n} S_{j,m} \, dy_j
\end{equation}
and
\begin{equation} \label{p:em}
  \begin{cases} P_1 = 0      \\
     P_m = \sum_{j < l}^{n} 
   \left( \sum_{k=1}^{m-1} \{ S_{j,k},S_{l, m-k} \} \right) \, dy_j \wedge dy_l 
          \quad\text{when} \ m\geq 2. 
   \end{cases}   
\end{equation}
where $\{ \cdot , \cdot \}$ denotes the Poisson bracket on $\bar{Z}$. Then
on every fibre $\bar{F}_b$, $S_m$ and $P_m$ 
satisfy the following equations
\begin{equation} \label{u:rec}
           d \, S_m|_{\bar{F}_b} = P_{m}|_{\bar{F}_b}.
\end{equation}
\end{prop}
\begin{proof}
We recall that the Poisson bracket on $T^*U / \Lambda$ can be written as
\[ \{ f, g \} = \sum_{k=1}^{n} 
           \partial_{y_k} f \, \partial_{b_k} g - 
                     \partial_{b_k}f \, \partial_{y_k}g .\]
Since the functions $S_{j,m}$ do not depend on $y_1$, one can easily see that 
the following holds
\[ \{ S_{j,k} \, b_1^k , \, S_{l,m} \, b_1^m \} = \{ S_{j,k} , \, S_{l,m} \} 
                                                         \, b_1^{k+m}, \] 
where the bracket on the right hand side reduces to the bracket on $\bar{Z}$.
Also we have
\[ \{ S_{j,k} \, b_1^k, \, o(b_1^N) \} = o(b_1^{N+k}). \]
Thus we have
\begin{eqnarray*}
    \{ u_j, u_l \} & = &
            \sum_{0 \leq m+k \leq N} \{ S_{j,k},\, S_{l,m} \} \, b_1^{m+k} +
                       o(b_1^{N}) \\
    & = &  \sum_{m = 0}^{N} 
             \left( \sum_{k=0}^{m} \{ S_{j,k},\, S_{l, m-k} \} \right) \, b_1^m +
                       o(b_1^N). 
\end{eqnarray*}
Therefore if $u_j$ and $u_l$ commute then we must have that for all $m \in \N$
\[ \sum_{k=0}^{m} \{ S_{j,k},\, S_{l, m-k} \} = 0, \]
or that
\begin{equation} \label{u:a}
    \{ S_{j,m}, S_{l,0} \} + \{ S_{j,0}, S_{l,m} \} = - 
                        \sum_{k=1}^{m-1} \{ S_{j,k},\, S_{l, m-k} \}.
\end{equation}
The condition that $u|_{Z} = \pi$ implies
\[ S_{j,0} = b_j. \]
Therefore
\[ \{ S_{j,m}, S_{l,0} \} = \{ S_{j,m}, b_l \} = \partial_{y_l} S_{j,m}. \]
We then see that (\ref{u:a}) becomes
\[ \partial_{y_l} S_{j,m} - \partial_{y_j} S_{l,m}  = - 
                        \sum_{k=1}^{m-1} \{ S_{j,k},\, S_{l, m-k} \}\]
which is exactly what we get by expanding (\ref{u:rec}).
\end{proof}
\begin{rem} Notice that (\ref{u:rec}) are a set of partial differential 
equations satisfied by the sequence $\{ S_m \}_{m \in \N}$. 
Moreover the definition of $P_m$ depends only on the $S_k$'s with $k \leq m-1$, 
therefore one may think of solving the equations recursively. From each 
solution $S_m$ of the $m$-th equation, we may determine another by adding to $S_m$ 
a fibrewise closed section of $\mathfrak{L}^{\ast}$.
\end{rem}

Now we  provide a method to construct and characterise 
sequences $\{ S_m \}_{m\in \N}$ of solutions to (\ref{u:rec}). Suppose 
$(V,u) \in \mathscr U_{\bar{Z}}$ and let $W \subseteq u(V)$ be a neighbourhood 
of $\Gamma$.
Let $r \in \R$ be a parameter. For
$b = (0, b_2, \ldots, b_n) \in \Gamma$, let $(r,b)$ denote the point
$(r,b_2, \ldots, b_n) \in \R^n$. Given $(r,b) \in W$, denote by 
$L_{r,b}$ the fibre $u^{-1}((r,b))$. For every fibre $F_b \subset Z$ of 
$\pi$, consider the symplectomorphism
\begin{equation} \label{fb:symp}
 (y_1, \ldots, y_n, \sum_{k=1}^{n} x_k dy_k) 
            \mapsto (x_1, b_2 + x_2, \ldots, b_n + x_n, y_1, \ldots, y_n),
\end{equation}
between a neighbourhood of the zero section of $T^{\ast}F_b$ and 
a neighbourhood of $F_b$ in $V$. 
If $W$ is sufficiently small, for every $(r,b) \in W$, the Lagrangian 
submanifold $L_{r,b}$ will be the image of the graph of a closed $1$-form 
on $F_b$. Due to the $S^1$ invariance of $u$ and the fact that $u_1=b_1$, this 1-form has to be of the type
\[ r dy_1 + \ell(r,b), \]
where $\ell(r,b)$ is the pull back to $F_b$ of a closed one form on $\bar{F}_b$.  
Denote by $\ell(r)$ the smooth one parameter family of sections of 
$\mathfrak{L}^{\ast}$ such that $\ell(r)|_{\bar{F}_b} = \ell(r,b)$.
The condition $u|_Z=\pi$ implies that $\ell (0,b)=0$. Furthermore, the $N$-th order Taylor series expansion of 
$\ell(r)$ in the parameter $r$ can be written as
\begin{equation} \label{l:tay}
 \ell(r) =\sum_{k=1}^{N} \ell_k \, r^k + o(r^N), 
\end{equation}
where the $\ell_k$'s are fibrewise closed sections of $\mathfrak{L}^{\ast}$. 
We can write
\begin{equation}\label{eq:ell_k}
\ell_k = \sum_{j=2}^{n} a_{j,k} \, dy_j.
\end{equation}
The following Lemma is rather technical but straightforward, 
thus its proof may be skipped on first reading.
\begin{lem} \label{ls:form} 
Given $(V,u) \in \mathscr U_{\bar{Z}}$, let $\{ S_m\}_{m\in\N}$ 
be the sequence (\ref{eq:S_m}) of sections of $\mathcal L^*$ encoding the Taylor coefficients of $u$ and let 
$\{ \ell_m \}_{m \in \N}$ be the sequence of fibrewise closed
sections of $\mathfrak{L}^{\ast}$ constructed from $u$ as 
above. Then for every $m \in \N$, there exist formulae
\begin{equation} \label{ls:rec}
 a_{j,m} = - S_{j,m} + R_{j,m}, 
\end{equation} 
where $R_{j,m}$ is an explicit polynomial expression depending on the 
$S_{l,k}$'s 
and their derivatives in the $b_i$'s up to order $m-1$ and with 
$0 \leq k \leq m-1$. In particular $R_{j,1} = 0$.
Thus the sequence $\{ \ell_m \}_{m \in \N}$ uniquely determines 
the sequence $\{ S_m \}_{m \in \N}$ recursively and viceversa. 
\end{lem}
\begin{proof} 
First of all let us write
\[ \ell(r) = r \, dy_1 + \sum_{j=2}^{n} a_j(r) dy_j. \]
Then by definition
\begin{equation} \label{agei:tay}
  a_j(r) = \sum_{k=1}^{N} a_{j,k} \, r^k + o(r^N). 
\end{equation}
The $a_j$'s are functions of $(r, b, y)$, with $(b,y) \in \bar{Z}$, satisfying
by construction
\begin{equation} \label{sistema1}
  \begin{cases}
        u_1(r, b_2 + a_2, \ldots, b_n + a_n, y) = r, \\
        u_j(r, b_2 + a_2, \ldots, b_n + a_n, y) = b_j \quad\text{for all} \ 
                                      j=2, \ldots, n. 
   \end{cases} 
\end{equation}
When $W$ is sufficiently small and $(r,b) \in W$, this system can be solved
using the implicit function theorem to determine the $a_j$'s uniquely. 
We will now use it to compute the $a_{j,m}$'s and determine the formulae 
(\ref{ls:rec}). 

\medskip
Let $j=2, \ldots, n$, then from the system and the conditions on $u$ 
we obtain
\[ a_j|_{r=0} = 0 \]
and
\begin{equation} \label{a:1stder}
  \partial_{b_1}u_j + \sum_{k=2}^{n} \partial_{b_k} u_j \,
                                       \partial_r a_k  = 0 .
\end{equation}
When evaluating at $r=0$, using $u|_{Z} = \pi$, we get
\[ \partial_{b_1}u_j|_{r=0} + \partial_r a_j|_{r=0}  = 0, \]
i.e. that
\begin{equation} \label{as:1} 
           a_{j,1} = - S_{j,1}. 
\end{equation}
\medskip
Now we do the second order terms. Derivating (\ref{a:1stder}) we obtain
\[ \partial_{b_1}^2 u_j + \sum_{k=2}^{n} \partial_{b_1} \partial_{b_k} u_j \,
           \partial_r a_k + \sum_{k,l=2}^{n} \partial_{b_l} \partial_{b_k} u_j \,
                      \partial_r a_l \, \partial_r a_k +
                         \sum_{k=2}^{n} \partial_{b_k} u_j \,
                                       \partial_r^2 a_k  = 0. \]
Evaluating at $r=0$ we get
\[ \partial_{b_1}^2 u_j|_{r=0} + 
  \left( \sum_{k=2}^{n} \partial_{b_1} \partial_{b_k} u_j \, \partial_r a_k
                                                   \right) |_{r=0} + 
                                                   \partial_r^2 a_j|_{r=0}  = 0,\]
i.e. we obtain
\begin{equation} \label{as:2}
  a_{j,2} = - S_{j,2} +
              \sum_{k=2}^{n} \partial_{b_k} S_{j,1} \, S_{k,1}. 
\end{equation}
So we have that (\ref{ls:rec}) holds for $m=2$, where
\[ R_{j,2} =  \sum_{k=2}^{n} \partial_{b_k} S_{j,1} \, S_{k,1}. \]
For the terms of order greater that two we refer the reader to the Appendix in \S \ref{appendix}.
\end{proof}

\begin{rem}
We point out that (\ref{as:1}) shows that the definition of $\ell_1$ given in
this section coincides with the first order invariant defined in the previous 
section.
\end{rem}

One good reason to work with the sequence $\{ \ell_k \}_{k \in \N}$
rather than with the sequence $\{ S_{m} \}_{m \in \N}$ is that we can easily 
prove the following
\begin{prop} \label{lk:exist} 
Given any sequence $\{ \ell_k \}_{k\in \N}$ of fibrewise closed 
sections of $\mathfrak{L}^{\ast}$, there exists a smooth
$1$-parameter family $\ell(r)$ of fibrewise closed sections of $\mathfrak{L}^{\ast}$ 
such that (\ref{l:tay}) holds for every $N \in \N$. 
\end{prop}

The proof of this is based on the following general:
\begin{lem}\label{lem:borel} For any sequence of $C^\infty$ functions $\{\alpha_k:\R^p\rightarrow\R\}$, there is a $C^\infty$ function $f:\R\times\R^p\rightarrow\R$, such that $\alpha_k(x)=\partial_r^kf(r,x)|_{r=0}$, for all $k\in\N$.
\end{lem}
A proof of this Lemma in the case when $\{\alpha_k\}$ is a sequence of real numbers is hinted in \cite{Rudin} Exercise 13, page 384. It is an exercise to show that the method proposed there can be adapted to the case when $\alpha_k$ depends smoothly on a parameter $x\in\R^p$.

\begin{proof}[Proof of Proposition \ref{lk:exist}]
Let us first prove the statement assuming that all the $\ell_k$'s are fibrewise 
exact, i.e. there exists a sequence of functions $\{ f_k \}_{k \in \N}$
on $\bar{Z}$ such that
\[  \ell_k|_{\bar{F}_b} = d \, f_k|_{\bar{F}_b}. \]
We have $\bar{Z} \cong \R^{n-1} \times T^{n-1}$, where $T^{n-1}$ is the 
$(n-1)$-torus. Let $\{ U_{\alpha}, \phi_{\alpha} \}_{\alpha \in J}$ be 
a partition of unity on $T^{n-1}$. Define 
\[ f_{k,\alpha} = \sqrt{\phi_{\alpha}} f_{k}. \]
We apply Lemma \ref{lem:borel}, for every $\alpha \in J$, to the sequence $\{ f_{k, \alpha} \}_{k \in \N}$
lifted to the covering $\R^{n-1}$ of $T^{n-1}$. So there exists a $C^\infty$ function
$f_{\alpha}=f_{\alpha}(r)$ such that
\[ f_{\alpha}(r) = \sum_{k=1}^{N} f_{k, \alpha} \, r^k + o(r^N), \]
for every $N \in \N$.
Let
\[ f(r) = \sum_{\alpha \in J} \sqrt{\phi_{\alpha}} \, f_{\alpha}(r).\]
Then $f(r)$ descends to a smooth $1$-parameter family of functions on $\bar{Z}$.
Moreover
\begin{eqnarray*}
  f(r) & = & \sum_{k=1}^{N}  \left( \sum_{\alpha \in J} \sqrt{\phi_{\alpha}} 
                                         f_{k, \alpha} \right) \, r^k + o(r^N) \\
       & = &  \sum_{k=1}^{N} \left( \sum_{\alpha \in J}
                                   \phi_{\alpha} f_{k} \right) \, r^k + o(r^N) \\
        & = &  \sum_{k=1}^{N} f_{k} \, r^k + o(r^N). \\
\end{eqnarray*}
If we let $\ell(r)$ be the $1$-parameter family of sections of $\mathfrak{L}^{\ast}$
such that
\[  \ell(r)|_{\bar{F}_b} = d \, f(r)|_{\bar{F}_b}, \]
then we clearly have
\[ \ell(r) =\sum_{k=1}^{N} \ell_k \, r^k + o(r^N). \]

We now do the general case. There certainly is a sequence $\{ l_k \}_{k \in \N }$
of fibrewise constant sections of $\mathfrak{L}^{\ast}$ such that for every 
$k \in \N$, $\ell_k - l_k$ is fibrewise exact. Since $l_k$ is fibrewise constant
we can write
\[ l_k = \sum_{j=2}^{n} q_{j,k} \, dy_j, \]
where the $q_{j,k}$'s are fibrewise constant functions. Invoking Lemma \ref{lem:borel}, for every 
$j =2, \ldots, n$, there exists a family of fibrewise constant functions
$q_j(r)$ such that
\[ q_j(r) =\sum_{k=1}^{N} q_{j,k} \, r^k + o(r^N). \]
Let 
\[ l(r) = \sum_{j=2}^{n} q_j(r) dy_j, \]
and let $\tilde{\ell}(r)$ be the fibrewise exact family of forms such that
\[ \tilde{\ell}(r) = \sum_{k=1}^{N} (\ell_k - l_k)   \, r^k + o(r^N), \]
which exists from the previous step.
Define
\[ \ell(r) = \tilde{\ell} (r) + l(r). \]
One easily checks that (\ref{l:tay}) holds.
\end{proof}

The following is an existence result

\begin{thm} \label{l to s}
Let $\{ \ell_m \}_{m \in \N}$ be a sequence of fibrewise closed
sections of $\mathfrak{L}^{\ast}$  and let $\{ S_m\}_{m \in \N}$ be the sequence 
of sections of $\mathfrak{L}^{\ast}$ obtained recursively from $\{ \ell_m \}_{m \in \N}$
using formulae (\ref{ls:rec}) in  Proposition~\ref{ls:form}, then there exists 
$(V,u) \in \mathscr U_{\bar{Z}}$ such that for every $N \in \N$
\[ u_j = \sum_{k=0}^{N} S_{j,k}  \, b_1^k + o(b_1^N). \]
\end{thm}

\begin{proof}
Following Proposition~\ref{lk:exist}, given the sequence $\{ \ell_m \}_{m \in \N}$, 
we can construct $\ell(r)$, a smooth $1$-parameter family of fibrewise closed 
sections of $\mathfrak{L}^{\ast}$ satisfying (\ref{l:tay}).
We show that $\ell(r)$ can be used to construct the pair $(V,u)$. 
In fact the process
is the inverse of the one which led us to the construction of a family
$\ell(r)$ from a fibration $u$. The construction is identical to the one in the proof 
of Theorem~\ref{broken:constr}.
Denote $\ell(r,b) = \ell(r)|_{\bar{F}_b}$ and write
\[ \ell(r,b) = \sum_{j=2}^{n} a_j(r,b) dy_j, \]
where the $a_j(r,b)$'s are functions depending on $y$ and they satisfy
\begin{equation} \label{aj:tay}
  a_j(r,b) = \sum_{k=1}^{N} a_{j,k}(b) \, r^k + o(r^N). 
\end{equation}  
Let $L_{r,b}$ be the Lagrangian submanifold of $T^*U/ \Lambda $ which
is the image of the closed one form $$r dy_1 + \ell(r,b)$$ under the symplectomorphism 
(\ref{fb:symp}).
When $W \subseteq U$ is sufficiently small and $(r,b) \in W$,
then the submanifolds $L_{r,b}$ are the fibres of a Lagrangian fibration 
$u: V \rightarrow \R^n$. We describe $u$ explicitly and show that $(V,u) \in \mathscr U_{\bar{Z}}$ . Given $(b^{\prime}, y^{\prime}) \in V$, 
there exists a unique $(r,b) \in W$ such that $L_{r,b} \subset V$ and
$(b^{\prime}, y^{\prime}) \in L_{r,b}$, in fact $(r,b)$ can be determined
as functions of $(b^{\prime}, y^{\prime})$ by solving the non linear system
\begin{equation} \label{sistema2}
     \begin{cases} r = b_1^{\prime} \\
                 b_j + a_j(r, b, y^{\prime}) = b_j^{\prime} \quad\text{when} 
                              \ j=2, \ldots, n    
   \end{cases}   
\end{equation}
using the implicit function theorem.
We define
\[ {u}_1(b^{\prime}, y^{\prime}) = b_1^{\prime} \]
and, when $j=2, \ldots, n$
\begin{equation} \label{ugeii}
         u_j(b^{\prime}, y^{\prime}) = b_j(b^{\prime}, y^{\prime}).
\end{equation} 
We claim that the coefficients of the Taylor series expansion of $u_j$ in 
$b_1^{\prime}$ are exactly the coefficients $S_{j,m}$ obtained from the sequence 
$\{ \ell_m \}_{m \in \N}$ through formulae (\ref{ls:rec}). 
In fact notice that, by construction of $u_j$, the functions 
$a_j$ satisfy
\[ u_j(r, b_2 + a_2(r, b, y^{\prime}), \ldots, b_n + a_n(r, b, y^{\prime})) = b_j, \]
i.e. they are obtained from $u_j$ as the unique solution to system
(\ref{sistema1}) and therefore the claim follows from the proof of
Lemma~\ref{ls:form}.
\end{proof}

\begin{cor} \label{lclos:s}
 Let $\{ \ell_m \}_{m \in \N}$ and $\{ S_m\}_{m \in \N}$ be 
sequences of sections of $\mathfrak{L}^{\ast}$ with coefficients $a_{j,k}$ and $S_{j,k}$, respectively, related by formulae (\ref{ls:rec}). Then all the $\ell_m$'s are fibrewise closed
if and only if the sequence $\{ S_m\}_{m \in \N}$ satisfies equations
(\ref{u:rec}).
\end{cor}
\begin{proof} If all the $\ell_m$'s are fibrewise closed, then Theorem~\ref{l to s}
shows that there is a Lagrangian fibration $u: V \rightarrow \R^n$
whose Taylor coefficients are given by the sequence $\{ S_m\}_{m \in \N}$.
Being $u$ Lagrangian, the claim follows from Proposition~\ref{u:tay eq}.

\medskip
Suppose now that $\{ S_m\}_{m \in \N}$ satisfies equations (\ref{u:rec}).
We prove the claim by induction. First of all notice that when $m=1$,
$S_1 = - \ell_1$ and equation (\ref{u:rec}) implies that $\ell_1$ is 
fibrewise closed. Now suppose we have proved that $\ell_m$ is fibrewise closed 
for all $m \leq N$. Consider the sequence $\{ \tilde{\ell}_m \}_{m \in \N}$, 
where $\tilde{\ell}_m = \ell_m$ when $m \leq N$ and $0$ otherwise. 
Using formulae (\ref{ls:rec}), we construct the associated sequence 
$\{ \tilde{S}_m\}_{m \in \N}$. Since all the $\tilde{\ell}_m$'s are fibrewise closed,
 from the first part of this Corollary it follows that $\{ \tilde{S}_m\}_{m \in \N}$
satisfies equations (\ref{u:rec}). Denote by $\tilde{P}_m$ the $2$-forms 
in (\ref{p:em}) constructed from $\{ \tilde{S}_m\}_{m \in \N}$.
Now notice that
\[ \tilde{S}_m = S_m, \]
when $m \leq N$ and 
\[ \tilde{P}_m = P_m \]
when $m \leq N+1$.  Moreover, if we denote by $\tilde{R}_{j,k}$ the 
expressions $R_{j,k}$ appearing in (\ref{ls:rec}) applied to 
$\{ \tilde{S}_m\}_{m \in \N}$, then
\[ \tilde{R}_{j,N+1} = R_{j,N+1}, \]
where the right hand side denotes the same expression obtained using 
$\{ S_m\}_{m \in \N}$. Therefore formula (\ref{ls:rec}) with $m = N+1$
and the fact that $\tilde{\ell}_{N+1} = 0$, implies
\begin{equation} \label{sn1}
    \tilde{S}_{j,N+1} = \tilde{R}_{j,N+1} = R_{j,N+1}. 
\end{equation}
Define the one form
\[ R_{N+1} = \sum_{j=2}^{n} R_{j, N+1} dy_j. \]
Clearly (\ref{sn1}) says that
\[ \tilde{S}_{N+1} = R_{N+1} \] 
and that equation (\ref{u:rec}) for $\{ \tilde{S}_m\}_{m \in \N}$ when $m = N+1$
becomes
\begin{equation} \label{d:er}
  d R_{N+1}|_{\bar{F}_b} = \tilde{P}_{N+1}|_{\bar{F}_b} = P_{N+1}|_{\bar{F}_b}. 
\end{equation}
Using (\ref{ls:rec}) for $\{ S_m\}_{m \in \N}$ when $m= N+1$, we obtain
\[ d \ell_{N+1}|_{\bar{F}_b} = - dS_{N+1}|_{\bar{F}_b} + d R_{N+1}|_{\bar{F}_b}. \]
Now substituting (\ref{d:er}) and using the fact that (\ref{u:rec}) holds for 
 $\{ S_m\}_{m \in \N}$ when $m = N+1$ we obtain
\[d \ell_{N+1}|_{\bar{F}_b}=  - dS_{N+1}|_{\bar{F}_b} + P_{N+1}|_{\bar{F}_b} = 0, \]
which completes the proof.
\end{proof}

Finally we have the most general existence result

\begin{thm} \label{hoi:exist}
Let $\{ S_m \}_{m \in \N}$ be a sequence of sections of  $\mathfrak{L}^{\ast}$,
satisfying (\ref{u:rec}), then there exists $(V,u) \in \mathscr U_{\bar{Z}}$ 
such that for every $N \in \N$
\[ u_j = \sum_{k=0}^{N} S_{j,k}  \, b_1^k + o(b_1^N). \]
\end{thm} 

\begin{proof} It is just a matter of applying the previous results. In fact 
given the sequence $\{ S_m \}_{m \in \N}$ satisfying (\ref{u:rec}), using 
Proposition~\ref{ls:form} we construct the sequence $\{ \ell_m \}_{m \in \N}$,
whose terms are all fibrewise closed thanks to Corollary~\ref{lclos:s}. 
Finally we apply Theorem~\ref{l to s} to obtain $u$.
\end{proof}

Define the following sets:
\[
\begin{array}{l}
\mathscr{L}_{\bar{Z}} = \{ \{ \ell_m \}_{m \in \N }\mid \ell_m \ \textrm{is
a}\ C^\infty,\ \textrm{fibrewise closed section of} \ \mathfrak{L}^{\ast} \};\\
\\
\mathscr{S}_{\bar{Z}}=\{ \{ S_m \}_{m \in \N }\mid S_m \ \text{is
a}\ C^\infty\ \textrm{section of} \ \mathfrak{L}^{\ast}\ \textrm{satisfying (\ref{u:rec})}\}.
\end{array}
\]
Clearly Proposition~\ref{u:tay eq} gives a map 
$T: \mathscr{U}_{\bar{Z}} \rightarrow \mathscr{S}_{\bar{Z}}$,
assigning to $(V,u)$ the Taylor coefficients of $u$. We summarise the previous 
results in the following:
\begin{thm}\label{thm:sequences} There is a one to one correspondence between
the sets $\mathscr{L}_{\bar{Z}}$, $\mathscr{S}_{\bar{Z}}$ and 
 $T( \mathscr{U}_{\bar{Z}})$.  In particular from every sequence  
$\ell \in \mathscr L_{\bar{Z}}$ we can construct a unique
element $S(\ell) \in \mathscr{S}_{\bar{Z}}$ and an element $(V,u) \in 
\mathscr{U}_{\bar{Z}}$ such that $T(V,u) = S(\ell)$.
\end{thm}

Given two elements $(V,u)$ and $(\tilde{V}, \tilde{u}) \in\mathscr U_{\bar{Z}}$, 
we can construct two stitched Lagrangian fibrations
$(Y,f_u)$ and $(\tilde{Y}, f_{\tilde{u}})$ as in (\ref{u:st}).  
We recall that $f_u$ and $f_{\tilde{u}}$
are equivalent up to a change of coordinates on the base 
if $f_{\tilde{u}} = \phi \circ f_u$, where $\phi:W\rightarrow \phi(W)\subseteq\R^n$ is an admissible change of coordinates on the base. If we write $\phi = (\phi_1, \ldots, \phi_n)$, then $\phi$ must satisfy $\phi_1 = b_1$ and $\phi|_{U^+} = \I$. 

\medskip
Similarly, we say that two sequences $\ell, \tilde{\ell} \in \mathscr{L}_{\bar{Z}}$ are equivalent up to 
a change of coordinates in the base if they define fibrations $f_u$ and $f_{\tilde{u}}$ respectively 
which are equivalent up to a change of coordinates in the base.
We now describe this equivalence relation in terms of a group action. Given a change 
of coordinate map $\phi$ on the base satisfying the above properties,
we can consider its Taylor expansion in $b_1$ from the left, i.e. where 
the coefficients are given by left derivatives. For each component 
$\phi_j$, $j = 2, \ldots, n$, it can be written as
\[ \phi_j(b_1, \ldots, b_n)|_{W \cap U^-} = b_j + 
              \sum_{k=1}^{N} \Phi_{j,k}(b_2, \ldots, b_n) b_1^k + o(b_1^N). \]
The left Taylor coefficients of $\phi$ thus define a sequence
$\{ \Phi_m \}_{m \in \N}$, where $\Phi_m: \Gamma \rightarrow \R^{n-1}$ is a $C^\infty$ map whose components are 
$\Phi_m = (\Phi_{2,m}, \ldots, \Phi_{n,m})$. 

\begin{lem}\label{lem:germ-adm} Given any sequence  $\{ \Phi_m \}_{m \in \N}$ of smooth maps
$\Phi_m: \Gamma \rightarrow \R^{n-1}$ with components 
$\Phi_m = (\Phi_{2,m}, \ldots, \Phi_{n,m})$ there exists an admissible
change of coordinate map $\phi=(\phi_1, \ldots, \phi_n) $ defined on some neighbourhood 
$W$ of $\Gamma$ such that $\phi_1 = b_1$, $\phi|_{U^+ \cap W} = \I$ and
\[ \phi_j(b_1, \ldots, b_n)|_{U^- \cap W} = b_j + 
              \sum_{k=1}^{N} \Phi_{j,k}(b_2, \ldots, b_n) b_1^k + o(b_1^N). \]
for all $N \in \N$.
\end{lem}

\begin{proof} It follows from Lemma \ref{lem:borel}.
\end{proof}

Define the following set
\[ \mathscr{D}_{\Gamma} = \{ \{ \Phi_{m} \}_{m \in \N} \, | \, 
                   \Phi_m \in C^{\infty}(\Gamma, \R^{n-1}) \}. \]
We say that two admissible change of coordinate maps $\phi$ and $\phi'$ are equivalent if their corresponding left Taylor coefficients define the same element in $\mathscr{D}_{\Gamma}$. We call $\mathscr{D}_{\Gamma}$ the set of \textbf{germs of admissible change of coordinates}.
Given a germ $\Phi\in\mathscr{D}_{\Gamma}$ we say that an admissible change of coordinates $\phi$ is a representative of $\Phi$, if $\phi$ satisfies Lemma \ref{lem:germ-adm}.

\medskip
Composition of germs of admissible maps induces a group structure on $\mathscr{D}_{\Gamma}$, i.e.
given $\Phi, \Phi^{\prime} \in \mathscr{D}_{\Gamma}$, we define 
$\Phi \cdot \Phi^{\prime}$ to be the germ of the map $\phi \circ \phi^{\prime}$, where $\phi$ and $\phi^{\prime}$ are representatives of $\Phi$ and $\Phi^{\prime}$ respectively. It is easy to see that this product on $\mathscr{D}_{\Gamma}$ does not depend on the choice of representatives.

\medskip
The group $\mathscr D_{\Gamma}$ acts on the set $\mathscr{L}_{\bar{Z}}$ as follows.
Given $ \ell\in\mathscr{L}_{\bar{Z}}$ and $\Phi\in\mathscr D_{\Gamma}$, we define $\Phi \cdot \ell$ to be the sequence $\tilde{\ell}\in\mathscr{L}_{\bar{Z}}$ associated to the Lagrangian fibration $\tilde{u} = \phi \circ u$, where $\phi$ is a representative of $\Phi$ and $u$ is a 
Lagrangian fibration obtained from $\ell$ via Theorem~\ref{l to s}.

\begin{lem} The above action is well-defined.
\end{lem}
\begin{proof} We need to show that the action does not depend on the choices made. The sequence  $\tilde{\ell}\in\mathscr{L}_{\bar{Z}}$ determines a unique sequence  $\{\tilde S_m\}_{m\in\N}\in\mathscr S_{\bar Z}$, where each $\tilde S_m$ is defined in terms of the Taylor coefficients of $\tilde u=\phi\circ u$. These coefficients, in turn, are expressed in terms of the Taylor coefficients of $\phi$ and $u$. If we take a different representative $\phi'$ of $\Phi$, clearly, the Taylor coefficients of $\phi'\circ u$ and $\phi\circ u$ coincide. Now let $S\in\mathscr S_{\bar Z}$ be the sequence determined by $\ell$. If $u'\in\mathscr U_{\bar Z}$ is a different realisation of $\ell$ then, by construction, $u'$ defines a sequence $S'\in\mathscr S_{\bar Z}$ such that $S=S'$. Therefore the Taylor coefficients of $\phi\circ u$ and $\phi\circ u'$ coincide.
\end{proof}

\begin{prop}
Let $f_u$ and $f_{\tilde{u}}$ be two stitched Lagrangian fibrations 
constructed as in (\ref{u:st}).  If locally $f_{\tilde{u}}= \phi \circ f_u$ for some admissible change of coordinates, then the sequences 
$\ell, \tilde{\ell} \in \mathscr{L}_{\bar{Z}}$ associated to $f_u$ and $f_{\tilde{u}}$
are in the same orbit of $\mathscr D_{\Gamma}$. Moreover $f_u$ is equivalent to a 
smooth fibration up to a change of coordinates on the base 
if and only if $\ell$ is a sequence of fibrewise constant sections of 
$\mathfrak{L}^{\ast}$.
\end{prop}
\begin{proof}
The first part of the statement is obvious. If $f_{\tilde{u}} = \phi \circ f_u$
is smooth, then $\tilde{\ell}$ is the zero sequence $0$. It is easy to verify
that $\ell = \Phi^{-1} \cdot 0$ is a sequence of fibrewise constant 
sections of $\mathfrak{L}^{\ast}$. Suppose, viceversa, that $\ell$ is a 
sequence of fibrewise constant sections. Consider the associated sequence 
$S \in \mathscr S_{\bar Z}$. The coefficients $S_{j,m}$ of each element
$S_m \in S$ can be regarded as functions on the base $\Gamma$, therefore
$S$ also defines a sequence $\Phi \in \mathscr D_{\Gamma}$ by setting
$\Phi_m = (S_{2,m}, \ldots, S_{n,m})$. Let $\phi$ be an admissible change
of coordinates representing $\Phi$. It is clear that $\phi^{-1} \circ f_u$ 
is smooth and that $\Phi^{-1} \cdot \ell = 0$. 
\end{proof}
In \S \ref{sec:normal} we will consider equivalences up to smooth fibre 
preserving symplectomorphism.

\section{The semiglobal classification}\label{sec:normal}
Let $(X,\omega )$ be a symplectic manifold and $f:X\rightarrow B$ be a stitched fibration as in Definition \ref{defi stitched}. Let $Z\subset X$ be the seam of $f$ and let $\Gamma:=f(Z)\subset B$. Since we are 
interested in a semiglobal classification, throughout this section we will consider
stitched Lagrangian fibrations satisfying the following assumption
\begin{ass} \label{st:ass}
The stitched Lagrangian fibration $f:X \rightarrow B$ satisfies
the following condition
\begin{enumerate}
\item the pair $(B, \Gamma)$ is diffeomorphic to the pair $(D^n, D^{n-1})$
      where $D^n \subset \R^n$ is a ball centred at the origin and 
      $D^{n-1} \subset D^n$ is the intersection of $D^n$ with an $n-1$ 
      dimensional subspace; 
\end{enumerate}
Also, the following data is specified
\begin{enumerate}
\setcounter{enumi}{1}
\item a basis $\gamma = (\gamma_1, \gamma_2, \ldots, \gamma_n)$ of $H_1(X, \Z)$
 so that $\gamma_1$ is represented by the orbit of the $S^1$ action;
\item a continuous section $\sigma$ of $f$ defined on a neighbourhood of 
$\Gamma$, such that $\sigma|_{f(X^+)}$ and $\sigma|_{f(X^-)}$ are restrictions 
of smooth maps on $B$ and the image of $\sigma$ is a smooth Lagrangian 
submanifold of $X$.\end{enumerate}
We denote a stitched Lagrangian fibration together with this data by 
$(X, B, f, \gamma, \sigma)$. 
\end{ass}
We will soon show that a section $\sigma$ as in $(3)$ always exists. 

\begin{defi}\label{def:symp-eq}
We say that two stitched fibrations 
$(X, B, f, \gamma, \sigma)$ and $(X', B', f', \gamma', \sigma')$, 
with seams $Z$ and $Z'$  respectively are 
\textbf{fibrewise symplectically equivalent} (or just equivalent) if 
there are neighbourhoods $ W \subseteq B$ of $\Gamma :=f(Z)$ and $W'\subseteq B'$ of $\Gamma':= f'(Z')$ 
and a commutative 
diagram:
\[ 
\begin{CD}
(f')^{-1}(W')  @> \Psi >> f^{-1}(W)\\
@Vf'VV   @VVfV \\
W'  @> \phi >> W 
\end{CD}
\]
where $\Psi$ is an $S^1$ equivariant $C^\infty$ symplectomorphism sending $Z'$ to $Z$ and $\phi$ is a $C^\infty$ diffeomorphism such that $\Psi\circ\sigma'=\sigma\circ\phi$ and
$\Psi_* \gamma' = \gamma$. The set of equivalence classes under this relation will be denoted by $\mathscr F$ and elements therein will be called \textbf{germs of stitched fibrations}.
\end{defi}

Now we show that any stitched fibration, satisfying 
Assumption~\ref{st:ass}, is fibrewise symplectically equivalent to a stitched fibration 
of the type $(Y,f_u)$ studied in \S \ref{sec:ho-inv}. 
Before doing this we need some preliminary results.

\medskip
Recall we can write a stitched fibration as:
\[
f=\begin{cases}
       		f^+ \quad\text{on} \ X^+;\\
                f^- \quad\text{on} \  X^-,
       \end{cases} 
\]
where $f^\pm$ is the restriction of a $C^\infty$, $S^1$ invariant map to $X^\pm$ whose fibres are Lagrangian when restricted to $X^\pm$. As pointed out in \S \ref{sec:def-ex}, the fibres of such a map are not a priori required to be Lagrangian beyond $X^\pm$. Nevertheless we have the following: 

\begin{prop}\label{prop: ext} Let $(X,B,f)$ be a stitched fibration with seam $Z\subset X$, satisfying condition (1) of Assumption~\ref{st:ass}. 
Then there are neighbourhoods $V\subseteq X$ of $Z$ and $W\subseteq B$ of $\Gamma:= f(Z)$ and a $C^\infty$, proper Lagrangian fibration $\tilde f^+:V\rightarrow W$ such that $\tilde f^+|_{X^+\cap V}=f^+|_{X^+\cap V}$. The same is true for $f^-$.
\end{prop}
\begin{proof}
Define $f_0:Z\rightarrow\Gamma$ to be $f_0=f|_{Z}$. Consider the reduced space $\bar Z=Z\slash S^1$ with its reduced symplectic form $\omega_{red}$.  On $\R\times S^1\times\bar Z$ we define the symplectic form:
\[
\omega_{red}+ds\wedge dt
\]
where $(t,s)$ are coordinates on $\R\times S^1$. From the coisotropic neighbourhood theorem (cf. \cite{Salamon} \S3.3) there exists a function $\epsilon:\Gamma\rightarrow\R_{> 0}$, a neighbourhood $V\subset X$ of $Z$ and a $S^1$-equivariant symplectomorphism between $V$ and
\begin{equation}\label{eq coiso}
\{(t,s,p)\in\R\times S^1\times\bar Z\mid -\epsilon(\bar f(p))<t<\epsilon(\bar f(p))\}.
\end{equation}
In particular, the projection onto $\R$ corresponds to the moment map $\mu$ on $V$. Now, on the set in (\ref{eq coiso}),  we can define an ``auxiliary" smooth Lagrangian fibration given by 
\[
\tilde\pi (t,s,p)=(t,f_0(s,p)).
\]
Fix a basis $\gamma$ of $H_1(V,\Z)\cong H_1(S^1\times\bar Z,\Z)$, 
satisfying condition $(2)$ of Assumption~\ref{st:ass} and a smooth Lagrangian 
section of $\tilde\pi$. The action-angle coordinates map $\Theta$ associated 
to $\tilde\pi$, with respect to $\gamma$ and $\sigma$, together with 
(\ref{eq coiso}), induces a $C^\infty$ symplectomorphism 
\begin{equation}\label{eq ident}
 \tilde V:=T^\ast U\slash\Lambda \cong V
\end{equation}
for some open neighbourhood $U$ of $0\in\R^n$ with coordinates $(b_1,\ldots ,b_n)$, which are the action coordinates of $\tilde\pi$. The pull back to 
$\tilde{V}$ of the $S^1$ action on $V$
is given by translations along $db_1$ and the corresponding moment map is $b_1$. Pulling back $f|_{V}$ to  $\tilde V$ via the latter identification we obtain a stitched fibration --with abuse of notation-- defined by:
\begin{equation}\label{eq:stitched-prime}
       f = \begin{cases}
       		u^+ \quad\text{on} \ \tilde V^+;\\
                u^- \quad\text{on} \  \tilde V^-,
       \end{cases} 
\end{equation}
where $u^\pm$ is the pull back of $f^\pm$. It follows that $u^+|_Z=u^-|_Z=\pi|_Z$.

\medskip
What we gained so far is an identification which allows us to view $f|_{V}$ as a stitched fibration on the smooth symplectic manifold $\tilde V$, where global canonical coordinates exist. Now we can use the results of \S 
\ref{sec:ho-inv} to show that $u^+$ (equivalently, $u^-$) can be extended as required. This can be done as follows. Since $u^+$ is the restriction of a $C^\infty$ map to $\tilde V^+$, all the  derivatives of its function components with respect to $b_1$ exist. Evaluating them at $b_1=0$ produces a unique sequence in $\mathscr S_{\bar Z}$ which in turn induces a unique sequence in $\mathscr L_{\bar Z}$ and a smooth Lagrangian fibration $(\tilde{V},w) \in \mathscr U_{\bar Z}$ 
(where eventually we restricted to a smaller $\tilde{V}$) whose Taylor coefficients in $b_1$ evaluated at $b_1=0$ coincide with those of $u^+$ (cf. Theorem \ref{thm:sequences}). In particular, this allows us to define
\begin{equation}
       \tilde u^+ = \begin{cases}
       		u^+ \quad\text{on} \ \tilde V^+;\\
                w \quad\text{on} \  \tilde V^-,
       \end{cases} 
\end{equation}
obtaining an element $( \tilde{V}, \tilde u^+) \in\mathscr U_{\bar Z}$, where
 $\tilde u^+$ extends $u^+$. Observe that different choices of $w$ induce different smooth extensions of $u^+$, however, all such choices are obtained starting from the same sequence in $\mathscr S_{\bar Z}$ determined by the derivatives of $u^+$. Finally, pulling back $\tilde u^+$ to $V$ under the 
identification (\ref{eq ident}), and perhaps shrinking $V$, gives us the 
required $\tilde f^+$. One can use the same arguments to find a suitable smooth extension of $f^-$.
\end{proof}

\begin{cor} A section $\sigma$ of $(X,B,f)$ satisfying condition (3) 
of Assumption~\ref{st:ass} exists.
\end{cor}
\begin{proof} Perhaps after an admissible change of coordinates on the 
base, a smooth Lagrangian section of $\tilde{f}^+$ is also a section of $f$.
\end{proof}

Let  $(b_1,\ldots b_n)$ be coordinates on $U\subseteq\R^n$ and let $\Lambda$ be the integral lattice inside $T^\ast U$ generated by $db_1,\ldots ,db_n$. 
Consider $T^\ast U\slash\Lambda$ with its standard symplectic structure and let $\pi: T^\ast U\slash\Lambda \rightarrow U$ be the standard projection. 
For convenience we change slightly the usual notation. Let 
$\gnor =\{ b_1=0\}\cap U$, $U^+:=\{b_1\geq 0\}\cap U$, 
$U^-:=\{b_1\leq 0\}\cap U$,  $\znor := \pi^{-1}(\Gamma)$ and 
$\zbnor := \znor / S^1$. 
We assume that $(U, \gnor)$ is diffeomorphic to the pair $(D^n, D^{n-1})$. 
Given  $(V, u) \in \mathscr U_{\zbnor}$, we can construct a stitched 
Lagrangian fibration $(Y, f_u)$ as in (\ref{u:st}).
The zero section of $\pi$ also defines a section of $f_u$,
after perhaps an admissible change of coordinates on the base. We denote
this section by $\sigma_0$. 
As a basis of $H_1(Y, \Z)$ we take the basis $(db_1, \ldots, db_2)$ of 
$\Lambda$. We denote it by $\gamma_0$. 
Then $(Y, f_u(Y), f_u, \sigma_0, \gamma_0)$ satisfies Assumption~\ref{st:ass}.

\begin{defi}
Let $F := (X,B,f, \sigma, \gamma)$ be a stitched Lagrangian fibration with 
seam $Z$, satisfying Assumption~\ref{st:ass}. A stitched fibration
$F_u := (Y, f_u(Y), f_u, \sigma_0, \gamma_0)$ of the type
above is a \textbf{normal form} of $F$ if $F_u$ and $F$ define the same germ 
of a stitched fibration, according
to Definition~\ref{def:symp-eq}.
\end{defi}

Observe that the above is a normalisation of a $T^n$-fibred neighbourhood of 
the seam of $F$. In this sense, $F_u$ is a \emph{semi-global} normal form.

If $F_u = (Y, f_u(Y), f_u, \sigma_0, \gamma_0)$ is a normal form of 
$F := (X,B,f, \sigma, \gamma)$ and $Z$ is the seam of $F$, then $\znor$ of 
$F_u$ is nothing else but $Z$ expressed in action angle coordinates, and thus it 
is a normalisation of $Z$.
Since $\sigma_0$ and $\gamma_0$ are chosen canonically, we will from now
on omit to specify them and just denote the normal form by $F_u= (Y,f_u)$.

\begin{prop}\label{prop:normalform}
Every stitched Lagrangian fibration $(X,B,f)$ satisfying 
(1) of Assumption~\ref{st:ass} has a section $\sigma$ and a basis $\gamma$ 
as in (2) and (3) of Assumption~\ref{st:ass} such that 
$(X,B,f, \sigma, \gamma)$ has a normal form 
$(Y, f_u)$ .
\end{prop}
\begin{proof}
From Proposition \ref{prop: ext} we can assume there exist open neighbourhoods $V\subset X$ of $Z$ and $W \subseteq B$ of $\Gamma$ and a proper 
smooth Lagrangian fibration $\tilde f^+ :V \rightarrow W$ extending $f^+$. Now, fixing a basis $\gamma$ of $H_1(V,\Z)$ as in (3) of Assumption~\ref{st:ass}
and a smooth Lagrangian section $\sigma$ of $\tilde f^+$,  we obtain a unique 
symplectomorphism
\[
\Theta^+ :T^\ast U\slash\Lambda\rightarrow V 
\] 
given by the action-angle coordinates associated to $\tilde f^+$. 
Then by defining $u$ to be the pull back of 
$\tilde{f}^-$ under $\Theta^+$ one readily sees that $f$ transforms into a 
fibration of the type $(Y, f_u)$.  
\end{proof}

\begin{defi}\label{def:seq-inv} Let $(X,B,f, \sigma, \gamma)$ be a stitched 
fibration with a normal form $(Y,  f_u)$. 
Let $\ell\in\mathscr{L}_{\zbnor}$ be the unique sequence determined by $u$. We denote $\inv (f_u):=(\zbnor,\ell)$ and we call it the 
\textbf{invariants} of  $(Y, f_u)$.  The invariants of 
 $F=(X,B,f, \sigma, \gamma)$ are defined to be $\inv (F):=\inv(f_u)$. 
\end{defi}

\begin{prop} Let $F = (X,B,f, \sigma, \gamma)$ and $F'=
(X',B',f', \sigma', \gamma')$ be stitched fibrations with 
normal forms $(Y, f_u)$ and $(Y',f_{u'})$ defining invariants 
$\inv (F)$ and $\inv(F')$ respectively. If $F$ and $F'$ are fibrewise 
symplectically equivalent, then $\inv(F)=\inv (F')$.
\end{prop}
\begin{proof} Assume there is a commutative diagram as in Definition \ref{def:symp-eq}. To keep notation simple, let us assume $W=B$ and $W'=B'$. 
We have the diagram with commutative squares:
\begin{equation}\label{eq:squares}
\begin{CD}
Y @<\Theta<< X @>\Psi>> X' @>\Theta'>>Y' \\
@Vf_uVV @VfVV @Vf'VV @VVf_{u'}V\\
U @<a<< B @>\phi>> B' @>a'>> U'.
\end{CD}
\end{equation}
Let us concentrate on the outermost square of (\ref{eq:squares}) and define $\tilde\Psi=\Theta'\circ\Psi\circ\Theta^{-1}$ and $\tilde\phi=a'\circ\phi\circ a^{-1}$. We claim that:

\begin{itemize}
\item[\textit{(i)}] $\inv (f_{\tilde\phi\circ u})=\inv (f_u)$; and 
\item[\textit{(ii)}] $\inv (f_{u'\circ\tilde\Psi})=\inv (f_{u'})$.
\end{itemize}

Since $\tilde\phi\circ f_u=f_{u'}\circ\tilde\Psi$, \textit{(i)} and \textit{(ii)} would imply that $\inv(F)=\inv(F')$. It is clear that $\zbnor$ 
and $\zbnor'$ must coincide. Observe that $\tilde{\Psi}$, restricted to $Y^+$,
is a symplectomorphism onto $(Y')^+$
which commutes with the projections $\pi$ and $\pi'$ on $T^*U$ and
$T^*U'$ and sends the zero section to the zero section. Therefore we must
have $\tilde{\Psi}|_{Y^+} = \tilde{\phi}^*$. 

To prove \textit{(i)} observe that
${\tilde{\phi}}^\ast|_{T^*U^+}$ must send the lattice $\Lambda'$ defining $Y'$ 
to the lattice $\Lambda$ defining $Y$. From this it follows that 
$\tilde\phi|_{U^+}$ is the identity map and the restriction of 
$\tilde\Psi$ to $Y^+$ is also the identity map. Then $\tilde\phi\circ u|_{V^+}=u|_{V^+}$. From this and the smoothness of $\tilde\phi\circ u$ it follows that the sequences in $\mathscr S_{\zbnor}$ defined by $\tilde\phi\circ u$ and $u$ coincide. Hence $\inv(f_u)=\inv(\tilde\phi\circ f_u)$.  Similarly, to prove \textit{(ii)}, observe that $u'\circ\tilde\Psi|_{V^+}=u'|_{V^+}$. Since $u'$ and 
$\tilde{\Psi}$ are smooth it follows that $\inv (f_{u'\circ\tilde\Psi})=\inv(f_{u'})$.
\end{proof}

\begin{cor} The definition of the invariants of $F = (X,B,f, \sigma, \gamma)$ 
is independent on the choice of normalisation.
\end{cor}
\begin{proof} 
Suppose we have two normalisations $(Y,f_u)$ and $(Y',f_{u'})$.
Clearly $\zbnor$ and $\zbnor'$ must coincide. We can also assume, w.l.o.g. that
$Y = Y'$. What may be different are the maps
$u$ and $u'$ such that $f_u$ and $f_{u'}$ are two different 
normalisations of $f$ induced from different extensions $\tilde f^+$ of $f^+$.
Consider the invariants $\inv(f_u)$ and $\inv(f_{u'})$, respectively. Since $f_u$ and $f_{u'}$ are symplectically equivalent via $\tilde\Psi=\Theta'\circ\Theta^{-1}$ and $\tilde\phi=a'\circ a^{-1}$, it follows that $\inv (f_u)=\inv(f_{u'})$.
\end{proof}

\begin{prop} \label{inv:sympl}
Let $F = (X,B,f, \sigma, \gamma)$ and $F'=
(X',B',f', \sigma', \gamma')$ be stitched Lagrangian fibrations 
satisfying Assumption~\ref{st:ass}. If $\inv (F)=\inv (F')$ then 
$F$ is fibrewise symplectically equivalent to $F'$ .
\end{prop}
\begin{proof}

Let $(Y, f_u)$ and $(Y',f_{u'})$ be normal forms of $F$ and $F'$, respectively.
 We can assume, w.l.o.g., $Y = Y'$. Let $S_u$ and $S_{u'}$ be the series in $\mathscr S_{\zbnor}$ defined by $u$ and $u'$ respectively. By assumption 
$S_u=S_{u'}$. This allows us to find Lagrangian fibrations $(\bar{V},\bar u), 
(\tilde{V},\tilde u),(\tilde V', \tilde u') \in\mathscr U_{\zbnor}$ such that 
\[
\tilde u =
\begin{cases} 
u\quad\textrm{on}\ \tilde V ^- \\
\bar u\quad\textrm{on}\ \tilde V^+
\end{cases}
\ \textrm{and}\quad
\tilde u' =
\begin{cases} 
u'\quad\textrm{on}\ (\tilde V')^- \\
\bar u\quad\textrm{on}\ (\tilde V')^+
\end{cases}
\]
where $S_{\bar u}=S_{\tilde u}=S_{\tilde u'}$. Now there is a neighbourhood $W$ 
of $\gnor$ and smooth symplectomorphisms $\Theta :T^\ast W\slash\Lambda\rightarrow \tilde V $ and $\Theta':T^\ast W\slash\Lambda\rightarrow \tilde V'$ which 
are the action-angle coordinate map of the fibrations $\tilde u$ and $\tilde u'$, respectively. Defining $\Psi=\Theta'\circ\Theta^{-1}$, it is clear that $\Psi|_{\tilde V^+}$ is the identity. Furthermore, when restricted to $\tilde V^-$, $\Psi$ sends the fibres of $\tilde u|_{\tilde V^-}=u|_{\tilde V^-}$ to the fibres of $\tilde u'|_{(\tilde V')^-}=u'|_{(\tilde V')^-}$. Therefore $\Psi$ is fibre preserving with respect to $f_u$ and $f_{u'}$. It follows that $f$ and $f'$ are symplectically equivalent.
\end{proof}

We summarise the previous Propositions in the following:

\begin{thm}\label{thm: grosso}  Let $F = (X,B,f, \sigma, \gamma)$ and $F'=
(X',B',f', \sigma', \gamma')$ be stitched Lagrangian fibrations 
satisfying Assumption~\ref{st:ass}, with invariants $\inv (F)$ and $\inv (F')$,
respectively. Then $F$ and $F'$ define the same germ if and only if 
$\inv (F)=\inv (F')$. In other words, the set of germs of stitched fibrations 
$\mathscr F$ is classified by the pairs $(\zbnor, \ell)$, where 
 $\ell \in \mathscr L_{\zbnor}$.
\end{thm}

The above provides a semi-global classification of stitched Lagrangian fibrations. In contrast to what happens for smooth Lagrangian submersions where no semi-global symplectic invariants exist, stitched fibrations in general do give rise to non trivial semi-global invariants.

We can now also state a more precise version of Theorem~\ref{broken:constr}: 

\begin{thm} \label{broken:constr2} \label{I added this theorem, 
for completeness. Check the statement, see if you agree!}
Let $(U, \Gamma)$ be a pair, where $U$ is an open neighbourhood 
of $0 \in \R^n$ and $\Gamma = U \cap \{ b_1 = 0 \}$. Assume $(U, \Gamma)$ 
is diffeomorphic to the pair $(D^n, D^{n-1})$.
Let $\bar{\Lambda} \subseteq T^{\ast} \Gamma$ be the lattice spanned by
$\{ db_2, \ldots, db_n \}$, and let
$\bar{Z} = T^{\ast} \Gamma / \bar{\Lambda}$, with projection
$\bar{\pi}: \bar{Z} \rightarrow \Gamma$ and bundle
$\mathfrak{L} = \ker \bar{\pi}_{\ast}$. Given integers $m_2, \ldots, m_n$
and a sequence $\ell = \{ \ell_k \}_{k \in \N} \in \mathscr L_{\bar{Z}}$
such that
\begin{equation} \label{int:cond2}
  \int_{db_j} \ell_1 = m_j \ \ \ \ \text{for all} \ \ j=2, \ldots, n,
\end{equation}
there exists a smooth symplectic manifold $(X, \omega)$ and a stitched
Lagrangian fibration $f: X \rightarrow U$ satisfying the following
properties:
\begin{itemize}
\item[\textit{(i)}] the coordinates $(b_1, \ldots, b_n)$ on $U$
          are action coordinates of $f$ with $\mu = f^{\ast}b_1$ the moment map
          of the $S^1$ action;
    \item[\textit{(ii)}] the periods $\{ db_1, \ldots, db_n \}$,
          restricted to $U^{\pm}$
          correspond to bases 
          $\gamma^{\pm} = \{ \gamma_1, \gamma_2^{\pm}, \ldots,
          \gamma_n^{\pm} \}$ of $H_1(X, \Z)$ satisfying $(i)$ and $(ii)$ of
          Corollary~\ref{broken:per};
    \item[\textit{(iii)}]  there is a Lagrangian section 
          $\sigma$ of $f$, such that $(\bar{Z}, \ell)$ are the invariants
          of $(X,f, U, \sigma, \gamma^+)$.
\end{itemize}
The fibration $(X,f, U)$ satisfying the above properties is unique up to 
fibre preserving symplectomorphism. 
\end{thm}
\begin{proof}
The construction of $(X, \omega)$ is like in the proof of 
Theorem~\ref{broken:constr}, i.e. 
\[ X = (X^+ \cup V^+) \cup_{\tilde{Q}} (X^- \cup V^-). \]
But now the map $u$, used to construct $\tilde{Q}$, is chosen so that 
the fibration $f_u: X^+ \cup V^+ \rightarrow \R^n$, defined by
\[ f_u =  \begin{cases} \pi^+ \quad\text{on} \ b_1 \geq 0 \\
                      u \quad\text{on} \ b_1 \leq 0.
        \end{cases}   \]
satisfies $\inv(f_u) = (\bar{Z}, \ell)$. Such a $u$ exists thanks to
Theorem~\ref{thm:sequences}. The fibration $f$ is again defined by
\[ f =  \begin{cases} \pi^+ \quad\text{on} \ X^+ \\
                      \pi^- \quad\text{on} \ X^-.
        \end{cases}   \]
It is clear that by construction $(X,f,U)$ satisfies $(i)-(iii)$. Notice that $\tilde{Q}$ matches the zero 
section of $\pi^+$ to the zero section of $\pi^-$. Therefore the
section $\sigma$ is just given by the zero section of $\pi^+$ on $U^+$ and 
by the zero section of $\pi^-$ on $U^-$. 

It is clear from the results proved in this Section (in particular from the 
existence of a normal form) that any stitched Lagrangian fibration 
$(X,f,U)$ satisfying $(i)-(iii)$ can be constructed in this way.

Uniqueness of $(X,f,U)$ is proved as follows. 
The only choice involved in the construction is the function $u$. 
Any other choice $u'$ must still satisfy $\inv(f_{u'})=(\bar{Z}, \ell)$.
Denote by $X$ and $X'$ the manifolds obtained from choices $u$ and $u'$
respectively. Let $\Psi: X \rightarrow X'$ be the map defined to be
the identity on $X^+$ and on $X^-$. One can see that $\Psi$ is well defined since the first order invariants of $f_{u}$ and $f_{u'}$ coincide. It is clearly a smooth symplectomorphism away from $Z$. We need to show that it is
smooth on $Z$. To see this we can use an argument similar to the one used in 
Proposition~\ref{inv:sympl}. If we think of $\Psi$ in the 
coordinates on $X^+ \cup V^+$,  $\Psi$ 
is a symplectomorphism sending the fibres of $f_u$ to the 
fibres of $f_{u'}$ and the zero section to the zero section.  
In a neighbourhood of $Z$ and in these coordinates, we can describe $\Psi$,  
as follows.
Since $\inv (f_u) = \inv (f_{u'})$, we can replace $u$ and $u'$ 
with $\tilde{u}$ and $\tilde{u}'$ as in Proposition~\ref{inv:sympl}.
Let $\Theta: T^*W/ \Lambda \rightarrow V^+$ and $\Theta': T^*W / \Lambda 
\rightarrow V^+$ be action angle coordinates of $\tilde{u}$ and $\tilde{u}'$ 
respectively, associated to the zero section and to the basis $\gamma^-$. 
Then, in these coordinates, $\Psi$ coincides with $\Theta' \circ \Theta^{-1}$.
It is therefore smooth. 
\end{proof}

\section{Stitched Lagrangian fibrations with monodromy}\label{sec. stitch w/mono}
We now study stitched Lagrangian fibrations defined over a non 
simply connected open set $U$. In this case it may be that the fibration 
has non-trivial monodromy. When the fibration is smooth, this monodromy 
is usually detected by the behaviour of the periods of the fibration expressed 
in terms of smooth coordinates on the base. In the case of stitched Lagrangian
fibrations there may not exist smooth coordinates on $U$, i.e. coordinates 
with respect to which the fibration is smooth. We will see how to detect 
monodromy from the behaviour of the first order invariant $\ell_1$. 
This will be done mainly through the discussion of examples.

In Example \ref{broken focus focus}, the fibration is topologically isomorphic
to a focus-focus fibration. The singular fibre is over $0 \in \R^2$.
Restricted to $X - f^{-1}(0)$, $f$ is a stitched Lagrangian fibration onto 
$U = \R^2 - \{ 0 \}$. We know that the locally constant presheaf on $U$
given by 
\[  W \mapsto H_1(f^{-1}(W), \Z) \]
has monodromy around $0$, i.e. the monodromy map
\[ \mathcal{M}_{b}: \pi_1(U) \rightarrow H_1(F_b, \Z) \]
at a fibre over $b \in U$ is non-trivial. In fact, if $e$ is a generator
of $\pi_1(U)$, $\mathcal{M}_{b}(e)$
is conjugate to the matrix
\[ \left( \begin{array}{cc} 1 & 1 \\
                     0 & 1 
    \end{array} \right). \]
We now look at a more general $2$-dimensional case.
\begin{ex}  \label{two:mon}
Let $U \subset \R^2$ be an open annulus in $\R^2$ centred at the origin.
As usual denote $U^+ = U \cap \{ b_1 \geq 0 \}$, $U^- = U \cap \{ b_1 \leq 0 \}$ 
and $\Gamma = U^+ \cap U^-$. This time $\Gamma$ is disconnected.
We let $\Gamma_{u} = \Gamma \cap \{ b_2 \geq 0 \}$ and 
$\Gamma_{d} = \Gamma \cap \{ b_2 \leq 0 \}$ be the upper and lower parts of 
$\Gamma$ respectively. 
Now let $f: X \rightarrow \R^2$ be a stitched Lagrangian fibration such that
$f(X) = U$. Observe that the seam $Z$ has two connected 
components: $Z_u = f^{-1}(\Gamma_u)$ and $Z_d = f^{-1}(\Gamma_d)$. 
Denote by $\bar{Z}_u$ and $\bar{Z}_d$ the respective $S^1$ quotients, i.e.
the connected components of $\bar{Z}$.
Given $b \in \Gamma_u$ and choosing a curve going anticlock-wise once around $0$ 
as generator $e \in \pi_1(U)$, suppose that with respect to a basis
 $\{ \gamma_1, \gamma_2 \}$ of $H_1(F_b, \Z)$ the monodromy is
\begin{equation} \label{monodr}
  \mathcal{M}_{b}(e) = \left( \begin{array}{cc} 1 & -m \\
                     0 & 1 
\end{array} \right),
\end{equation}
for some integer $m \neq 0$. In this case we must have that
$\gamma_1$ is represented by the orbits of the $S^1$ action. 
As usual let $X^{\pm} = f^{-1}(U^{\pm})$. 
Since $U - \Gamma_{d}$ is contractible we can
think of $\{ \gamma_1, \gamma_2 \}$ as a basis of $H_1(f^{-1}(U-\Gamma_{d}), \Z)$.
Consider the diagrams:

\[
\xymatrix{
 &  H_1(X^{+},\Z) \ar[dr]  \\
 H_1(f^{-1}(U-\Gamma_d),\Z) \ar[ur] \ar[rr]^{j_+} & & H_1(f^{-1}(U-\Gamma_u),\Z)
}
\]
or
\[
\xymatrix{
H_1(f^{-1}(U-\Gamma_d),\Z) \ar[dr] \ar[rr]^{j_-} & & H_1(f^{-1}(U-\Gamma_u),\Z) \\
&  H_1(X^{-},\Z) \ar[ur]
}
\]
induced by inclusions and restrictions.  The map $j_+$ identifies 
$\{ \gamma_1, \gamma_2 \}$ with a basis $\{ \gamma_1, \gamma_2^+ \}$ 
of $H_1(f^{-1}(U-\Gamma_{u}), \Z)$, whereas $j_-$ with a basis
$\{ \gamma_1, \gamma_2^- \}$. 
Notice that monodromy is given by $j_+^{-1}\circ j_-$.  Therefore we must have $\gamma_2^+ = m \gamma_1 + \gamma_2^-$.
Hence $\{ \gamma_1, \gamma_2^+ \}$ and $\{ \gamma_1, \gamma_2^- \}$ satisfy 
conditions $(i)$ and $(ii)$ of Corollary~\ref{broken:per}.  
Applying Lemma~\ref{broken:action} to $f$ restricted to  $f^{-1}(U-\Gamma_u)$ 
we can consider the action coordinates map $\alpha$ constructed
by taking action coordinates with respect
to $\{ \gamma_1, \gamma_2^+ \}$ on $U^+$ and with respect to 
$\{ \gamma_1, \gamma_2^- \}$ on $U^-$.
Denote by $(b_1^d, b_2^d)$ such coordinates. Similarly on $U-\Gamma_{d}$ we 
can consider action angle coordinates with respect to the basis
$\{ \gamma_1, \gamma_2 \}$. Denote by $(b_1^u, b_2^u)$ these coordinates.
In particular we can identify
\[ \bar{Z}_d = T^{\ast} \Gamma_d \, /  \, \langle db_2^d \rangle_{\Z} \]
and
\[ \bar{Z}_u= T^{\ast} \Gamma_u \, / \, \langle db_2^u \rangle_{\Z} \]
With respect to this choice of coordinates we can construct
the first order invariants $\ell_1^u$ and $\ell_1^d$
of $f$ on $\bar{Z}_u$ and $\bar{Z}_d$ respectively, then
by applying  Remark~\ref{rem:action} we obtain
\[ \int_{db_2^u} \ell_1^u = 0 \ \ \text{and} \ \ \int_{db_2^d} \ell_1^d = m. \]
This tells us that monodromy can be read from a jump in cohomology class
of the first order invariant associated to action coordinates.
\end{ex}

Using the methods of Theorem~\ref{broken:constr}
we can also construct stitched Lagrangian fibrations with prescribed monodromy and
and invariants. In fact we have

\begin{thm} Let $U \subset \R^2$ be an annulus as above with coordinates $(b_1, b_2)$. Let
$\bar{Z}_d = T^{\ast} \Gamma_d \, /  \, \langle db_2 \rangle_{\Z}$ 
and
$\bar{Z}_u= T^{\ast} \Gamma_u \, / \, \langle db_2 \rangle_{\Z}$
with projections $\bar{\pi}^d$ and $\bar{\pi}^u$ and bundles
 $\mathfrak{L}_d = \ker\bar{\pi}^d_{\ast}$ and $\mathfrak{L}_u = \ker\bar{\pi}^u_{\ast}$
 respectively. Given an integer $m$ and sequences 
$\ell^d = \{ \ell_k^d \}_{k \in \N} \in \mathscr L_{\bar{Z}_d}$ and 
$\ell^u = \{ \ell_k^u \}_{k \in \N} \in \mathscr L_{\bar{Z}_u}$ such that
\[ \int_{db_2} \ell_1^u = 0 \ \ \text{and} \ \ \int_{db_2} \ell_1^d = m, \]
there exists a smooth symplectic manifold $(X, \omega)$ and a stitched 
Lagrangian fibration $f: X \rightarrow U$ having monodromy (\ref{monodr}) with
respect to some basis $\gamma = \{ \gamma_1, \gamma_2 \}$ of $H_{1}(f^{-1}(U- \Gamma_{d}), \Z)$
and satisfying the following properties:
  \begin{itemize}
    \item[\textit{(i)}] the coordinates $(b_1, b_2)$
          are action coordinates of $f$ with moment map $f^{\ast}b_1$; 
    \item[\textit{(ii)}] the periods $\{ db_1, db_2 \}$, restricted to $U^{\pm}$ 
          correspond to the basis $\{ \gamma_1, \gamma_2  \}$; 
    \item[\textit{(iii)}]  there is a Lagrangian section 
          $\sigma$ of $f$, such that $(\bar{Z}_u, \ell^u)$ and $(\bar{Z}_d, \, \ell^d)$ are the invariants
          of $(f^{-1}(U-\Gamma_d),\, f, \, U-\Gamma_d, \,  \sigma, \, \gamma)$ and
          $(f^{-1}(U- \Gamma_u),\, f, \, U-\Gamma_u, \, \sigma, \, j_+(\gamma))$ respectively.
\end{itemize} 
The fibration $(X,f, U)$ satisfying the above properties is unique up to 
fibre preserving symplectomorphism. 
\end{thm}
\begin{proof} We let $\Lambda_+$ and $\Lambda_-$ be the lattices generated
by $db_1$ and $db_2$ in  $T^{\ast}U^+$ and $T^{\ast}U^-$ respectively.
Define $X^{\pm} = T^{\ast}U^{\pm} / \Lambda_{\pm}$, 
$Z^{\pm}_{u} = (\pi^{\pm})^{-1}(\Gamma_u)$ and 
$Z^{\pm}_{d} = (\pi^{\pm})^{-1}(\Gamma_d)$. Then, using $\ell_1^u$ and
$\ell_1^d$, we construct maps
\[ Q_u: Z^{-}_{u} \rightarrow Z^+_{u} \]
and 
\[ Q_d: Z^{-}_{d} \rightarrow Z^+_{d} \]
like in Theorem \ref{broken:constr}. We use these maps to glue 
$X^+$ and $X^-$ topologically along their boundary and thus form
$X$. For the smooth and symplectic gluing we follow the same method as in 
Theorem \ref{broken:constr2}, where higher order invariants are used. 
From the discussion of Example~\ref{two:mon} it follows that the fibration has the 
prescribed monodromy. Uniqueness is proved like in Theorem \ref{broken:constr2}.
\end{proof}

We now discuss a three dimensional example.

\begin{ex} \label{amoeb:mon}
In $\R^3$ consider the three-valent graph 
\[ \Delta = \{(0,0,- t), \ t \geq 0 \} \cup \{ (0,- t,0), \  t \geq 0 \} 
             \cup \{ (0,t,t), \ t\geq 0 \} \]
and let $D$ be a tubular neighbourhood of $\Delta$. Take $U = \R^3 - D$
and assume we have a stitched Lagrangian fibration $f: X \rightarrow \R^3$
such that $U=f(X)$. The seam is $Z= f^{-1}( \{ b_1 = 0 \} \cap U)$. 
Again we let 
$U^+ = U \cap \{ b_1 \geq 0 \}$, $U^- = U \cap \{ b_1 \leq 0 \}$ 
and $\Gamma = U^+ \cap U^-$. Also let $X^{\pm} = f^{-1}(U^{\pm})$. 
This time $\Gamma$ (and thus $Z$) has three connected components
\begin{eqnarray*}
  \Gamma_c & = & \{ ( 0, t,s), \ t,s < 0 \} \cap U, \\
  \Gamma_d & = & \{ (0, t,s), \ t> 0, s < t \} \cap U, \\
  \Gamma_e & = & \{ (0, t,s), \ s> 0, t < s \} \cap U.
\end{eqnarray*}
Also denote by $Z_c$, $Z_d$ and $Z_e$ the corresponding connected
components of $Z$ and by $\bar{Z}_c$, $\bar{Z}_d$ and $\bar{Z}_e$
their $S^1$ quotients.

Fix $b \in \Gamma_c$ and suppose that there is a basis 
$\{ \gamma_1, \gamma_2, \gamma_3 \}$ of $H_1( F_b, \Z)$
and generators $e_0, e_1, e_2$ of $\pi_1(U)$, satisfying $e_0 e_1 e_2 = 1$, with
respect to which the monodromy transformations are
\begin{equation}
\mathcal{M}_b(e_1) =  T_1  =  \left( \begin{array}{ccc}
                 1 & 0 & -m_1 \\
                 0 & 1  & 0 \\
                 0 & 0  & 1 
              \end{array} \right), \ \ \ 
\mathcal{M}_b(e_2) =  T_2  = \left( \begin{array}{ccc}
                 1 & -m_2 &  0 \\
                 0 & 1  & 0 \\
                 0 & 0  & 1 
              \end{array} \right). \label{t12}
\end{equation}
and ${M}_b(e_0) = T_0 = T_1^{-1} T_2^{-1}$, for non zero integers $m_1$ and
$m_2$.  We have that $\gamma_1$ is represented by the orbits of the
$S^1$ action, since it is the only monodromy invariant cycle. 
Now, since $U - ( \Gamma_d \cup \Gamma_e)$ is
contractible, $\{ \gamma_1, \gamma_2, \gamma_3 \}$ is
a basis of $H_1( f^{-1}(U - ( \Gamma_d \cup \Gamma_e)), \Z)$.
Consider the diagrams:
\[
\xymatrix{
 &  H_1(X^{+},\Z) \ar[dr]  \\
 H_1(f^{-1}(U-(\Gamma_d\cup\Gamma_e)),\Z) \ar[ur] \ar[rr]^{j_+} & & H_1(f^{-1}(U-(\Gamma_c\cup\Gamma_d)),\Z)
}
\]
or
\[
\xymatrix{
H_1(f^{-1}(U-(\Gamma_d\cup\Gamma_e)),\Z) \ar[dr] \ar[rr]^{j_-} & & H_1(f^{-1}(U-(\Gamma_c\cup\Gamma_d)),\Z) \\
&  H_1(X^{-},\Z) \ar[ur]
}
\]
induced by inclusions and restrictions. The map $j_+$ identifies
$\{ \gamma_1, \gamma_2, \gamma_3 \}$ with a basis  of 
$H_1(f^{-1}(U- (\Gamma_c \cup \Gamma_d)), \Z) $, which 
we call $\{ \gamma_1, \gamma_2^+, \gamma_3^+ \}$, while $j_-$ 
identifies it with another basis, which we call 
$\{ \gamma_1, \gamma_2^-, \gamma_3^- \}$. 
Notice that the monodromy map $\mathcal{M}_b(e_1) = j_+^{-1}\circ j_-$. We must have 
\begin{equation} \label{tre:mon} 
 \begin{cases}
  \gamma_2^+ = \gamma_2^-, \\
  \gamma_3^+ = m_1 \gamma_1 + \gamma_3^- .
 \end{cases}
\end{equation}
Therefore $\{ \gamma_1, \gamma_2^+, \gamma_3^+ \}$ and
$\{ \gamma_1, \gamma_2^-, \gamma_3^- \}$ satisfy conditions
$(i)$ and $(ii)$ of Corollary~\ref{broken:per}. 
Applying Lemma~\ref{broken:action} to $f$ restricted to 
$f^{-1}(U- (\Gamma_c \cup \Gamma_d))$, we can consider the action coordinates 
map $\alpha$ on $U- (\Gamma_c \cup \Gamma_d) $ constructed by taking action 
coordinates with respect to $\{ \gamma_1, \gamma_2^+ ,\gamma_3^+ \}$ on $U^+$ 
and with respect to $\{ \gamma_1, \gamma_2^-, \gamma_3^- \} \}$ on $U^-$.
Let us denote these coordinates by $(b_1^e, b_2^e, b_3^e)$. Similarly we can 
consider action coordinates on $U - (\Gamma_d \cup \Gamma_e)$
with respect to the basis $\{ \gamma_1, \gamma_2, \gamma_3 \}$  of 
$H_1(f^{-1}(U- (\Gamma_d \cup \Gamma_e)), \Z) $. We denote them
by $(b_1^c, b_2^c, b_3^c)$.
We have the identifications
\[ \bar{Z}_e = T^{\ast} \Gamma_e \, /  \, \langle db_2^e, db_3^e  \rangle_{\Z} \]
and
\[ \bar{Z}_c= T^{\ast} \Gamma_c \, / \, \langle db_2^c, db_3^c \rangle_{\Z}. \]

With respect to these coordinates we can compute the first order invariants 
$\ell_1^e$ and $\ell_1^c$ on $\bar{Z}_e$ and $\bar{Z}_c$ respectively.
From Remark~\ref{rem:action} and identities (\ref{tre:mon}) 
applied to $\ell_1^c$ and 
$\ell_1^e$ we obtain
\begin{equation*}
  \int_{db_2^c} \ell_1^c  
=  \int_{db_3^c} \ell_1^c  = 0
\end{equation*}
and
\begin{equation*}
  \int_{db_2^e} \ell_1^e =  0 \ \ \text{and} \ \ 
  \int_{db_3^e} \ell_1^e  =  m_1.
\end{equation*}
Similarly we construct the first order invariant $\ell_1^d$ on $\bar{Z}_d$. 
It will satisfy
\begin{equation*}
  \int_{db_2^d} \ell_1^d  =  m_2 \ \ \text{and} \ \ 
  \int_{db_3^d} \ell_1^d =  0.    
\end{equation*}
Again, monodromy is understood in terms of the difference in the
cohomology class of the first order invariant. Example \ref{ex amoebous fibr} 
is a special case of this situation, where $m_1 = m_2 = 1$. 
\end{ex}

Again, one can produce stitched Lagrangian fibrations of the type described
in this example with the gluing method Theorem~\ref{broken:constr}. In fact
we can prove
\begin{thm} Let $U \subset \R^3$, $\Gamma_c$, $\Gamma_d$ and $\Gamma_e$
be as in Example~\ref{amoeb:mon} and let $(b_1, b_2, b_3)$ be coordinates
on $U$. Define
$\bar{Z}_c = T^{\ast} \Gamma_c \, /  \, \langle db_2, db_3 \rangle_{\Z}$,
$\bar{Z}_d= T^{\ast} \Gamma_d \, / \, \langle db_2, db_3 \rangle_{\Z}$
and $\bar{Z}_e= T^{\ast} \Gamma_e \, / \, \langle db_2, db_3 \rangle_{\Z}$
with projections $\bar{\pi}^c$, $\bar{\pi}^d$, $\bar{\pi}^e$ and bundles
$\mathfrak{L}_c = \ker\bar{\pi}^c_{\ast}$, $\mathfrak{L}_d = \ker\bar{\pi}^d_{\ast}$,
$\mathfrak{L}_e = \ker\bar{\pi}^e_{\ast}$. 
Suppose we are given integers $m_1$, $m_2$ and sequences
$\ell^c = \{ \ell_k^c \}_{k \in \N} \in \mathscr L_{\bar{Z}_c}$, 
$\ell^d = \{ \ell_k^d \}_{k \in \N} \in \mathscr L_{\bar{Z}_d}$ and 
$\ell^e = \{ \ell_k^e \}_{k \in \N} \in \mathscr L_{\bar{Z}_e}$ satisfying 
\begin{eqnarray*}
  \int_{db_2} \ell_1^c & = & \int_{db_3} \ell_1^c  = 0, \\
  \int_{db_2} \ell_1^e & = &  0 \ \ \text{and} \ \ 
  \int_{db_3} \ell_1^e  =  m_1, \\
  \int_{db_2} \ell_1^d & = & m_2 \ \ \text{and} \ \ 
  \int_{db_3} \ell_1^d =  0.     
\end{eqnarray*}
Then there exists a smooth symplectic manifold $(X, \omega)$ and a stitched 
Lagrangian fibration $f: X \rightarrow U$ having the same monodromy
of Example~\ref{amoeb:mon} with respect to some basis 
$\gamma = \{ \gamma_1, \gamma_2, \gamma_3 \}$ of 
$H_1( f^{-1}(U - ( \Gamma_d \cup \Gamma_e)), \Z)$   and satisfying the following 
properties:
\newcounter{mon3}
  \begin{list}{(\roman{mon3})}{\usecounter{mon3} 
                            \setlength{\parsep}{0cm} 
                            \setlength{\topsep}{\itemsep}
                            \setlength{\leftmargin}{.5cm}}
    \item[\textit{(i)}] the coordinates $(b_1, b_2, b_3)$
          are action coordinates of $f$ with moment map $f^{\ast}b_1$; 
    \item[\textit{(ii)}] the periods $\{ db_1, db_2, db_3 \}$, restricted to $U^{\pm}$ 
          correspond to the basis $\gamma$; 
    \item[\textit{(iii)}]  there is a Lagrangian section 
          $\sigma$ of $f$, such that $(\bar{Z}_c, \ell^c)$, $(\bar{Z}_d, \, \ell^d)$ and 
          $(\bar{Z}_e, \, \ell^e)$ are the invariants
   of $(f^{-1}(U-( \Gamma_d \cup \Gamma_e)),\, f, \, U-( \Gamma_d \cup \Gamma_e), \,  \sigma, \, \gamma)$,
   $(f^{-1}(U-( \Gamma_c \cup \Gamma_e)),\, f, \, U-( \Gamma_c \cup \Gamma_e), \,  \sigma, \, j_+(\gamma))$
   and  $(f^{-1}(U-( \Gamma_c \cup \Gamma_d)),\, f, 
                    \, U-( \Gamma_c \cup \Gamma_d), \,  \sigma, \, j_+(\gamma))$ 
         respectively.
 \end{list} 
The fibration $(X,f, U)$ satisfying the above properties is unique up to 
fibre preserving symplectomorphism. 
\end{thm}

We omit the proof which is simply a repetition of the usual gluing method
from Theorems~\ref{broken:constr} and \ref{broken:constr2}.

\section{More examples?}
In this section we would like to propose a conjectural construction
generalising the one, described in \cite{CB-M-torino},  which led us to
Example~\ref{ex amoebous fibr}. In \cite{Gui-Stern-bi}, Guillemin and
Sternberg make the following observation. Let $N=n+m$, with $n,m$ positive
integers. Consider $\C^{N+1}$ with its standard symplectic structure, then
$S^1$ acts on it, in a Hamiltonian way, via the action given by,
\begin{equation} \label{simple:act}
 \theta:(z_1, \ldots, z_{N+1} )
  \mapsto (e^{i\theta} z_1, e^{-i\theta} z_2, \ldots, e^{-i\theta} z_{n+1},
             z_{n+2}, \ldots, z_{N+1})
\end{equation}
with moment map
$$\mu = \frac{|z_1|^2 - |z_2|^2 -  \ldots - |z_{n+1}|^2}{2}.$$
The action is singular along $\Sigma = \{ z_1= \ldots =z_{n+1}=0 \}$,
which can be identified with $\C^m$.
The observation is that for any $\epsilon \in \R_{\geq 0}$
the reduced spaces $(M_{\epsilon}, \omega_{r}(\epsilon))$ can be identified
with $(\C^N, \omega_{\C^N})$ with standard symplectic form,
(this includes the case of the critical value $\epsilon =0$). While when
$\epsilon \in \R_{<0}$, $(M_{\epsilon}, \omega_r(\epsilon))$ can be
identified with the $\epsilon$-blow up of $(\C^N, \omega_{\C^N})$ along
the symplectic submanifold $\Sigma$.

\medskip
The $\epsilon$-blow up can be described as follows. Let $L$ be the
total space of the tautological line bundle on $\PP^{n-1}$. The incidence
relation gives $L$ as
\[ L = \{ (v, l ) \in \C^{n} \times \PP^{n-1} \, | \, v \in l \}. \]
There are two natural projections: $\pi: L \rightarrow \PP^{n-1}$,
which is the bundle projection, and $\beta: L \rightarrow \C^{n}$ which is
the
blow-up map. The latter is a biholomorphism onto $\C^{n} - \{0\}$ once
the zero section is removed from $L$. Let $\omega_{FS}$ be the standard
Fubini-Study symplectic form on $\PP^{n-1}$. The $\epsilon$-blow up of
$\C^n$ at $0$ is $L$ together with the symplectic form given by
\[ \omega_{\epsilon} = \beta^{\ast} \omega_{\C^n} + \epsilon \,
                                      \pi^{\ast} \omega_{FS}. \]
The $\epsilon$-blow up of $\C^N$ along $\Sigma = \C^m$ can be
identified with $L \times \C^m$ with symplectic form
$\omega_{\epsilon} + \omega_{\C^m}$.

\medskip
In the case $n=1$ the blow-up is topologically (and holomorphically)
trivial, i.e. blowing up does not do anything. In fact one can also
show, by following Guillemin and Sternberg's argument, that the
reduced spaces can all be identified with $(\C^{m+1}, \omega_{\C^{m+1}})$
for all values of $\epsilon$.
This identification can also be explained as follows. Consider the
map $\gamma$ given in (\ref{eq. g}) and define the map
$p: \C^{m+2} \rightarrow \C^{m+1}$
given by
\begin{equation} \label{pws:p}
     p: (z_1, z_2, z_3, \ldots, z_{m+2}) \mapsto (\gamma(z_1, z_2), z_3,
 \ldots, z_{m+2}).
\end{equation}
Restricted to $\mu^{-1}(\epsilon)$, this map can be regarded as the quotient
map $\mu^{-1}(\epsilon) \rightarrow M_\epsilon$. It can be shown that the
reduced symplectic form with respect to this map is precisely
$\omega_{\C^{m+1}}$.

\medskip
Example~\ref{ex amoebous fibr} comes from this construction in the case
$m=n=1$. In fact the fibration $f$ is of the type $\Log \circ \Phi \circ p$,
where $\Phi$ is a symplectomorphism of $\C^2$ and
$\Log:(\C^{\ast})^2 \rightarrow \R^2$ is the map
$(v_1, v_2) \mapsto (\log|v_1|, \log |v_2|)$. The fact that $f$ is not
smooth is due to the non-smoothness of $p$, i.e. the reduced spaces are
not identified with $\C^2$ in a smooth way.

\medskip
We think that it may be possible to generalize this construction. The idea
is to use another result of Guillemin and Sternberg proved in the same
paper. The result is as follows. Let $\bar{X}$ be a compact $2N$ dimensional
symplectic manifold with symplectic form $\omega$ and $2m$-dimensional
symplectic submanifold $Y$. Consider now a principal $S^1$ bundle
$p_0: P \rightarrow \bar{X}$ with a connection one form $\alpha$. Given an
interval $I=(-\epsilon, \epsilon)$, Guillemin and Sternberg \cite{Gui-Stern-bi}\S 12 construct a
$2(N+1)$ symplectic manifold $X$ with the following properties.
\begin{enumerate}
\item There exists a Hamiltonian $S^1$ action on $X$ with
proper, surjective moment map $\mu: X \rightarrow I$.
\item For positive $t \in I$, $\mu^{-1}(t)$ is equivalent, as an $S^1$
bundle,
to $P$ and the reduced symplectic space $(X_t, \omega_r(t))$ is
symplectomorphic to $(\bar{X}, \omega)$.
\item The only critical value of $\Phi$ is $t=0$. If
$\Sigma := \Crit(\mu) \subset \mu^{-1}(0)$, i.e. the set of critical
points of
$\mu$, then $\Sigma$ is a smooth symplectic, $2m$ dimensional
submanifold of $X$ and the $S^1$ action is locally modelled on
(\ref{simple:act}) (in this case $0$ is also called a simple critical value).
If $X_0$ denotes the reduced symplectic space at $0$, with reduced symplectic
form $\omega_{r}(0)$ and quotient map $\pi_0: \mu^{-1}(0) \rightarrow X_0$,
then the triple $(X_0, \pi_0(\Sigma), \omega_r(0))$ can be identified with
$(\bar{X}, Y, \omega)$.
\item When $t \in I$ is negative, then the reduced space
$(X_t, \omega_r(t))$ can be identified with the blow-up $\tilde{X}$ of
$\bar{X}$ along $Y$ with symplectic form $\omega_{Y,t} + \beta^* t d\alpha$,
where $\omega_{Y,t}$ is the $t$-blow-up form along $Y$ on $\tilde{X}$ and
$\beta: \tilde{X} \rightarrow \bar{X}$ is the blow down map.
\end{enumerate}

We are interested in Guillemin-Sternberg's construction in the case $N =
m+1$,
i.e. in the case $Y$ is a codimension $2$ symplectic manifold. For simplicity
we also assume that $P= \bar{X} \times S^1$ and $\alpha=0$.
We can make the following observations.
\begin{itemize}
\item[(a)]
      Topologically $\tilde{X}$ is equivalent to $\bar{X}$, but
symplectically
      $(\tilde{X}, \omega_{Y,t})$ and $(\bar{X}, \omega)$ differ since the
      latter one has less area (blowing up removes the area of a small
tubular
      neighbourhood of $Y$).
\item[(b)]
      Consider the quotient $p: X \rightarrow X / S^1$, then $X / S^1$ can
      be identified with $\bar{X} \times I$. If we restrict $p$ to
      $X - \Sigma$ then it becomes an $S^1$ bundle onto
      $(\bar{X} \times I) - (Y \times \{0 \})$.
      Let $c_1$ be the first Chern class of this bundle. If $S$ is
      a small $2$-sphere centred at the origin in a fibre of the normal
      bundle of $Y  \times \{ 0 \}$ inside $(\bar{X} \times I)$, then
      $c_1(S) = 1$.
\end{itemize}

As we saw in the beginning of this section, in the non-compact case
$(\bar{X}, \omega) = (\C^{m+1}, \omega_{\C^{m+1}})$ and $Y = \C^m$,
the observation in $(a)$ was not true, in the sense that the identification
could be made also symplectically. This is because, although blowing up
locally reduces area, in this non-compact case the area is infinite so
it does not constitute a symplectic invariant. So the idea is to try to
generalize Guillemin and Sternberg's construction to other non-compact
cases. One interesting situation is if we take $(\bar{X}, \omega)$
with  $\bar{X} = (\C^*)^N$  and
\[ \omega = \sum_{k=1}^{N} \frac{dz_k \wedge d\bar{z}_k}{|z_k|^2}. \]
As symplectic submanifold $Y$ we can take some smooth algebraic hypersurface.

We think it may be possible to generalize Guillemin and Sternberg's
construction to this case. The hypothesis of compactness was made in order to
be able to use the coisotropic embedding theorem in symplectic topology, but
this theorem holds also in non-compact situations.
The question is whether the reduced spaces can all be identified with
$((\C^*)^N, \omega)$. Since the space is non-compact, area
is not an obstruction.

\begin{figure}[!ht]
	\centering
	\includegraphics{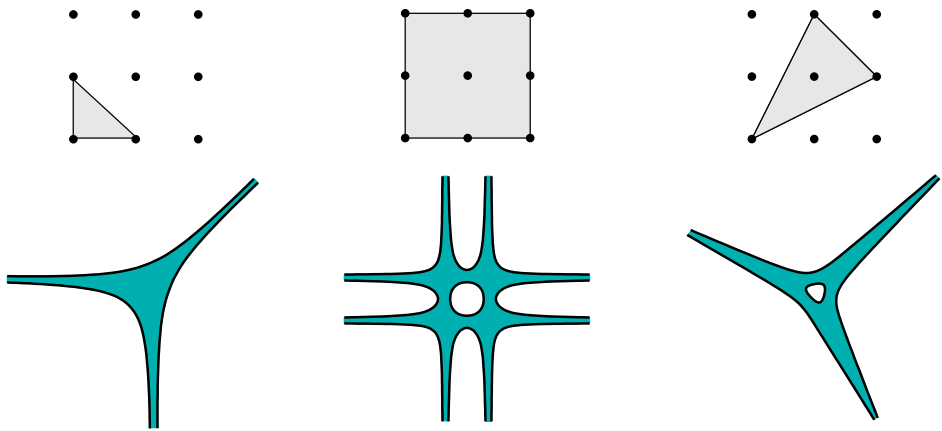}
	\caption{Amoebas with their respective Newton polygons.}
	\label{fig: general_amoeba}
\end{figure}

\medskip
Why would such a construction be useful? We could use it to construct
interesting examples of piecewise smooth Lagrangian fibrations with singular
fibres. In fact suppose the conjectured symplectic manifold $X$ exists
with the above properties and such that all reduced spaces can be identified
with $((\C^*)^N, \omega)$. Then, on $X$ we could define a piecewise smooth
Lagrangian fibration as follows.
On $(\C^*)^N \times I$ define the $T^N$ fibration given by
\[ F : (z_1, \ldots, z_N, t) \rightarrow ( \log|z_1|, \ldots, \log|z_N|,
t). \]
Clearly $F_t = F|_{(\C^*)^N \times \{t \}}$ is Lagrangian.
Now suppose there exists a map $p: X \rightarrow (\C^*)^N \times I$,
equivalent to the quotient $X \rightarrow X/S^1$ and with respect to which
the reduced spaces are all $((\C^*)^N, \omega)$. Presumably this map would be
locally modelled on (\ref{pws:p}), in particular it would fail to be
smooth on
$\mu^{-1}(0)$. The piecewise smooth Lagrangian fibration would be
\begin{equation}
f = F \circ p.
\end{equation}
We expect $f$ to be a stitched Lagrangian fibration when restricted to
$X - f^{-1}(\Delta)$.
The interesting aspect of this map is the structure of the singular fibres.
In fact its discriminant locus is $\Delta = F(Y \times \{ 0 \})$, which is
$\Log(Y) \times \{0 \}$. Images of algebraic hypersurfaces of
$(\C^*)^N$  by $\Log$ are called amoebas and they have shapes of the
type pictured in Figure~\ref{fig: general_amoeba}

The topological property, discussed in the observation (b), of the
bundle $p: X - \Sigma \rightarrow (\bar{X} \times I) - (Y \times \{0\})$,
ensures that the fibration $f$, restricted to $X-f^{-1}(\Delta)$ has
non-trivial monodromy. In fact one can find examples where monodromy
would be of the types discussed in (\ref{amoeb:mon}).
These examples, and the calculation of monodromy, generalize the construction
in \cite{TMS} of the negative fibre, also called the fibre of type $(2,1)$,
where a circle bundle with the topological property $(b)$ is used.

\medskip
In a work in progress \cite{CB-M} the authors use the piecewise smooth Lagrangian fibration in Example~\ref{ex amoebous fibr} as one of the building blocks for the construction of Lagrangian 
fibrations of $6$-dimensional compact Calabi-Yau manifolds.
One of the ideas involved is that the invariants we have defined for stitched
Lagrangian fibrations can be used to perturb the fibration in Example~\ref{ex amoebous fibr} away
from the singular fibres in order to glue it to other pieces of fibration.
In fact the sequence $\ell = \{ \ell_k \}_{k \in \N}$ of fibrewise
closed sections of $\mathcal L^*$ on $\bar{Z}$ can be easily perturbed,
for example by multiplying each element by cut-off functions on the base $\Gamma$ or
by summing to each element other fibrewise closed section and so on.

We believe that the more general construction proposed in this section is interesting because, if it can be carried through, 
then these Lagrangian fibrations could be used as building blocks of more general Lagrangian fibrations of 
compact symplectic manifolds.

\section{Appendix to Lemma \ref{ls:form}}\label{appendix}

We give here a proof of Lemma \ref{ls:form} for all $m\in\N$. Recall we can write 
\begin{equation}\label{aj:tay_again}
  a_j(r) = \sum_{k=1}^{N} a_{j,k} \, r^k + o(r^N). 
\end{equation}
The $a_j$'s are functions of $(r, b, y)$, with $(b,y) \in \bar{Z}$, satisfying
\[
  \begin{cases}
        u_1(r, b_2 + a_2, \ldots, b_n + a_n, y) = r, \\
        u_j(r, b_2 + a_2, \ldots, b_n + a_n, y) = b_j \quad\text{for all} \ 
                                      j=2, \ldots, n. 
   \end{cases} 
\]
When $W$ is sufficiently small and $(r,b) \in W$, the functions  $a_{j,m}$'s can be uniquely determined using the implicit function theorem. We will now use it to compute the $a_{j,m}$'s and obtain formulae (\ref{ls:rec}). We can rewrite the second equation of the above system by applying 
\begin{equation} \label{ugei:tay_again}
 u_j = \sum_{k=0}^{N} S_{j,k} b_1^k + o(b_1^N). 
\end{equation}
We obtain
\[ b_j + a_j + \sum_{k=1}^{N} S_{j,k}( b_2 + a_2, \ldots, b_n + a_n, y) r^k + o(r^N) =
                     b_j \]
which implies
\begin{equation} \label{aj:bla}
 a_j + \sum_{k=1}^{N} S_{j,k}( b_2 + a_2, \ldots, b_n + a_n, y) r^k + o(r^N) =
                     0. 
\end{equation}
To express everything as a power series in $r$ we use the Taylor
expansion up to a certain order $N^{\prime}$ of the $S_{j,k}$'s, which in 
the multi-index notation is given by:
\[ S_{j,k}( b_2 + a_2, \ldots, b_n + a_n, y) = \sum_{l=0}^{N^{\prime}} \sum_{|I| = l}
        C_I \, \partial^l_{I} S_{j,k}(b_2, \ldots, b_n) 
                             \, a_{2}^{i_2} \cdot \ldots \cdot a_{n}^{i_n} + 
\ldots, \]
where $I=(i_2, \ldots, i_n)$ is a multi-index and the $C_I$'s are suitable 
constants. 

\medskip
Let us introduce the following notation. For every multi-index $I = (i_2, \ldots, i_n)$,
let us define the following set
\[ \mathcal{H}_{I} = \{ (H_2, \ldots, H_n) \, | \, H_k \in (\Z_{>0})^{i_k} \ 
\text{if} \ i_k \geq 1 \ \text{and} \ H_k = 0 \in \Z
                                     \ \text{if} \ i_k = 0 \}. \]
When $i_k \geq 1$, we also write $H_k = (h_{k,1}, \ldots, h_{k,i_{k}})$.
 For every $m \in \N$, we denote
\[ \mathcal{H}_{I, m} = \left\{ (H_2, \ldots, H_n) \in \mathcal{H}_{I} \, | 
                           \, \sum_{i_k \neq 0} \sum_{j=1}^{i_k} h_{k,j} = m \right\}. \]
Clearly if $|I| = 0$ and $m \geq 1$ or if $0 \leq m < |I|$ 
then $\mathcal{H}_{I, m}$ is empty. 
When $i_k \neq  0$ for all $k=2, \ldots, n$, substituting (\ref{agei:tay}) we compute 
that
\[ a_{1}^{i_1} \cdot \ldots \cdot a_{n}^{i_n} = 
     \sum_{m=1}^{N^{\prime}} \left( \sum_{H \in \mathcal{H}_{I, m}} 
              a_{2,h_{2,1}} \cdot \ldots \cdot  a_{2,h_{2,i_2}} \cdot \ldots \cdot
                       a_{n,h_{n,1}} \cdot \ldots \cdot  a_{2,h_{n,i_n}} \right) r^m 
                        + o(r^{N^{\prime}}). \]

Let us introduce another bit of notation. When $|I| \neq 0$, 
for all $H \in \mathcal{H}_I$, let
\[ A_{H} = \prod_{i_{k} \neq 0} \prod_{j=1}^{i_k} a_{k,h_{k,j}}. \]
When $|I|=0$, the only element in $\mathcal{H}_I$ is $0 \in \Z^n$, so we set
\[ A_0 = 1. \]
Thus for all multi-indices $I$, we have
\[ a_{1}^{i_1} \cdot \ldots \cdot a_{n}^{i_n} = 
     \sum_{m=0}^{N^{\prime}} \left( \sum_{H \in \mathcal{H}_{I, m}} 
              A_H  \right) r^m 
                        + o(r^{N^{\prime}}). \]
Therefore $S_{j,k}( b_2 + a_2, \ldots, b_n + a_n, y)$ written as a power
series in $r$ becomes
 \[ S_{j,k}( b_2 + a_2, \ldots, b_n + a_n, y) = \sum_{m=0}^{N^{\prime}} 
       \left( \sum_{|I| \leq m} \sum_{H \in \mathcal{H}_{I, m}} C_I \,
           \partial^{|I|}_{I} S_{j,k}(b,y) \, A_H \right)   r^m 
                        + o(r^{N^{\prime}}) . \]
Substituting this into (\ref{aj:bla}) we obtain
\[ a_j + \sum_{l=1}^{N} \left( \sum_{m=0}^{l-1} \sum_{|I| \leq m} 
             \sum_{H \in \mathcal{H}_{I, m}} C_I \,
           \partial^{|I|}_{I} S_{j,l-m}(b,y) \, A_H \right) r^l + o(r^N) = 0. \]
Substituting also (\ref{aj:tay_again}) we have
\[ \sum_{l=1}^{N}  \left( a_{j,l} +  \sum_{m=0}^{l-1} \sum_{|I| \leq m} 
             \sum_{H \in \mathcal{H}_{I, m}} C_I \,
           \partial^{|I|}_{I} S_{j,l-m}(b,y) \, A_H \right) r^l + o(r^N) = 0. \]
Therefore, for every $l \in \Z_{>0}$, we have
 \[ a_{j,l} = -  \sum_{m=0}^{l-1} \sum_{|I| \leq m} 
             \sum_{H \in \mathcal{H}_{I, m}} C_I \,
           \partial^{|I|}_{I} S_{j,l-m}(b,y) \, A_H. \]
When $l=1$, this becomes
\[  a_{j,1} = - S_{j,1}, \] 
when $l \geq 2$ it can also be written as
 \[ a_{j,l} = - S_{j,l} -  \sum_{m=1}^{l-1} \sum_{|I| \leq m} 
             \sum_{H \in \mathcal{H}_{I, m}} C_I \,
           \partial^{|I|}_{I} S_{j,l-m} \, A_H. \]
Now notice that when $1 \leq m \leq l-1$ and $H \in \mathcal{H}_{I, m}$,
then $A_H$ only depends on the $a_{j,k}$'s with $1 \leq k \leq l-1$. 
Therefore if we define
\[ R_{j,l} =  -  \sum_{m=1}^{l-1} \sum_{|I| \leq m} 
             \sum_{H \in \mathcal{H}_{I, m}} C_I \,
           \partial^{|I|}_{I} S_{j,l-m}(b,y) \, A_H, \]
when $l \geq 2$ and $R_{j,1} = 0$, then (\ref{ls:rec}) holds with
$R_{j,m}$ satisfying the required properties.

\bibliographystyle{plain}

\begin{thebibliography}{10}

\bibitem{Arnold}
V.I. Arnold.
\newblock {\em Mathematical Methods of Classical Mechanics}.
\newblock Springer Verlag, 2nd. edition, 1989.

\bibitem{CB-M}
R.~Castano-Bernard and D.~Matessi.
\newblock {Lagrangian 3-torus fibrations over singular affine manifolds}.
\newblock In progress.

\bibitem{CB-M-torino}
R.~Castano-Bernard and D.~Matessi.
\newblock {Some piece-wise smooth Lagrangian fibrations}.
\newblock {\em Rend. Sem. Mat. Univ. Politec. Torino \textbf{63}}, pages
  223--253, 2005.
\newblock Electronic: http://seminariomatematico.dm.unito.it or 
       in preprint version: http://www.mis.mpg.de/preprints/2005/prepr2005\_47.html.

\bibitem{Dui}
J.J. Duistermaat.
\newblock On global action-angle coordinates.
\newblock {\em Comm.Pure Appl. Math. \textbf{6}}, pages 678--706, 1980.

\bibitem{TMS}
M.~Gross.
\newblock {Topological Mirror Symmetry}.
\newblock {\em Invent. Math. \textbf{144}}, pages 75--137, 2001.
\newblock arXiv: math.AG/9909015.

\bibitem{G-Siebert}
M.~Gross and B.~Siebert.
\newblock {Affine manifolds, Log structures and Mirror Symmetry}.
\newblock {\em Turkish J. Math. \textbf{27}}, 2003.
\newblock arXiv: math.AG/0211094.

\bibitem{G-Siebert2003}
M.~Gross and B.~Siebert.
\newblock {Mirror Symmetry via Logarithmic degeneration data I}.
\newblock {\em J. Differential Geom. \textbf{72}}, pages 169--338, 2006.
\newblock arXiv: math.AG/0309070.

\bibitem{G-Wilson2}
M.~Gross and P.M.H. Wilson.
\newblock {Large Complex Structure Limits of $K3$ Surfaces}.
\newblock {\em J. Differential Geom. \textbf{55}}, pages 475--546, 2000.
\newblock arXiv: math.DG/0008018.

\bibitem{Gui-Stern-bi}
V.~Guillemin and S.~Sternberg.
\newblock Birational equivalence in the symplectic category.
\newblock {\em Invent. Math. \textbf{97}}, pages 485--522, 1989.

\bibitem{Hitchin}
H.~Hitchin.
\newblock {The moduli space of special Lagrangian submanifolds. Dedicated to
  Ennio De Giorgi}.
\newblock {\em Ann. Scuola Norm. Sup. Pisa Cl. Sci. (4) \textbf{25}}, pages
  503--515, 1997.
\newblock arXiv: dg-ga/9711002.

\bibitem{Joyce-SYZ}
D.~Joyce.
\newblock {Singularities of special Lagrangian fibrations and the SYZ
  conjecture}.
\newblock {\em Comm. Anal. Geom. \textbf{11}}, pages 859--907, 2003.
\newblock arXiv: math.DG/0011179.

\bibitem{Kontsevich-Soibelman}
M.~Kontsevich and Y.~Soibelman.
\newblock Homological mirror symmetry and torus fibrations.
\newblock In {\em Symplectic geometry and mirror symmetry (Seoul, 2000)}, pages
  203--263. World Sci. Publishing, River Edge, NJ, 2001.
\newblock arXiv: math.DG/0011041.

\bibitem{Salamon}
D.~McDuff and D.~Salamon.
\newblock {\em Introduction to symplectic topology}.
\newblock Oxford Mathematical Monographs, Clarendon, 1998.

\bibitem{Ruan}
W-D. Ruan.
\newblock {Lagrangian Torus Fibrations and Mirror Symmetry of Calabi-Yau
  Manifolds}.
\newblock In {\em Symplectic geometry and mirror symmetry (Seoul, 2000)}, pages
  385--427. World Sci. Publishing, River Edge, NJ, 2001.
\newblock arXiv: math.DG/0104010.

\bibitem{Rudin}
W.~Rudin.
\newblock {\em {Real and Complex Analysis}}.
\newblock McGraw-Hill Book Co., New York, third edition, 1987.

\bibitem{SYZ}
E.~Strominger, S-T. Yau, and E.~Zaslow.
\newblock {Mirror symmetry is $T$-duality}.
\newblock {\em Nucl. Phys. B \textbf{479}}, pages 243--259, 1996.
\newblock arXiv: hep-th/9606040.

\end{thebibliography}

\vspace{2cm}

\begin{flushleft}

Ricardo~CASTA\~NO-BERNARD \\
Max-Planck-Institut f\"ur Mathematik\\
Vivatsgasse 7, \\
D-53111, Bonn, Germany\\
e-mail: \texttt{castano@mpim-bonn.mpg.de}\\ \
\\ \ 
\  
\\ 
Diego~MATESSI\\
Dipartimento
 di Scienze e Tecnologie Avanzate\\
Universit\`{a} del Piemonte Orientale\\
Via Bellini 25/G\\
 I-15100 Alessandria, Italy\\
e-mail: \texttt{matessi@unipmn.it}\\

\end{flushleft}

\end{document}